%% file: Dias_enumeration_etds_final.tex
\documentclass{etds}

\usepackage{color}

\newtheorem{theo}{Theorem}[section]
\newtheorem{prop}{Proposition}[section]

\newcommand{\restrictto}[1]{\,\rule[-2.5mm]{0.125mm}{5.5mm}_{\rule[-1.0mm]{0mm}{
4mm}\, #1}}

\begin{document}
\ETDS{1}{13}{XX}{2011}

\runningheads{K.\ Dias}{Enumerating combinatorial classes of vector fields}

\title{Enumerating combinatorial classes of the
complex polynomial vector fields in $\mathbb{C}$}

\author{KEALEY DIAS}

\address{Mathematisches Seminar, Christian-Albrechts Universit{\"a}t zu Kiel, Ludewig-Meyn-Str. 4,\ 24098 Kiel,
Germany\\
Current address: Department of Mathematics and Computer Science,\\ Bronx Community College, 2155 University Ave., Bronx, New York 10453, USA\\
\email{kealey.dias@gmail.com}}

\recd{$5$ November $2010$}

\begin{abstract}
 In order to understand the parameter space $\Xi_d$ of monic and centered
complex polynomial vector fields in $\mathbb{C}$ of degree $d$, decomposed by the
combinatorial classes of the vector fields, it is interesting to know the
number of loci in parameter space consisting of vector fields with the same
combinatorial data (corresponding to topological classification with fixed separatrices at infinity).  \par
This paper answers questions posed by Adam L. Epstein and Tan Lei about the total
number of combinatorial classes and the
number of combinatorial classes corresponding to loci of a
specific (real) dimension $q$ in parameter space, for fixed degree $d$; these numbers are denoted by $c_d$ and $c_{d,q}$ respectively. These results are extensions of a result by Douady, Estrada, and Sentenac, which shows that the number of
combinatorial classes of the structurally stable complex polynomial vector fields in
$\mathbb{C}$ of degree $d$ is the Catalan number $C_{d-1}$. \par
We show that enumerating the combinatorial classes is equivalent to a so-called
\emph{bracketing problem}. Then we analyze the generating functions and find closed-form expressions for $c_d$ and $c_{d,q}$, and we furthermore make an asymptotic analysis of these sequences for $d$ tending to
$\infty$.\par
These results are also applicable to special classes of Abelian differentials, quadratic differentials with double poles, and singular holomorphic foliations of the plane.
\end{abstract}
%
\begin{figure}[b]%
\flushleft
\rule[0mm]{54.0mm}{0.15mm}\\ \hspace*{5.0mm}%
{\footnotesize  2000 Mathematics Subject Classification: 37F75, 05A15, 05A16.\par
 }
\end{figure}
\section{Introduction}
The space $\Xi_d$ of monic and centered single-variable complex polynomial vector fields of
degree $d$ is parameterized by its $d-1$ complex coefficients. The space $\Xi_d \simeq \mathbb{C}^{d-1}$ is decomposed into combinatorial classes where the
vector fields within each class have the same \emph{combinatorial data set} (to be defined). Each
of these combinatorial classes is a connected manifold with well-defined (real)
dimension $q$, which is the dimension of the combinatorial class as a subspace in $\Xi_d$. \par
Counting the combinatorial classes will help us to better
understand $\Xi_d$, and the techniques utilized may also prove valuable to those who are
considering other enumerative problems in dynamical systems. \par
In this paper, we study the properties of the sequences $c_d$ (the total number of
combinatorial classes in $\Xi_d$) and $c_{d,q}$ (the number of combinatorial classes of
dimension $q$ in $\Xi_d$). These will be referred to as the \emph{simplified
problem} and the
\emph{complete problem} respectively.
We first show in Section \ref{bracketingsetupsection} that these problems are equivalent
to so-called
\emph{bracketing problems}, where one counts the number of pairings of parentheses
satisfying desired conditions \cite{Com}. Then for each problem, a recursion
equation and
implicit expressions for the algebraic generating functions are calculated (Section
\ref{genfctsection}),
closed-form expressions for $c_d$ and $c_{d,q}$ are calculated (Section \ref{closedformssection}), and
asymptotic questions are
considered (Section \ref{asymptoticbehaviorsection}). In particular, we compute the
asymptotic growth of the
sequence $c_d$ and prove that the discrete probability distribution
$\pi_{d,q}=\frac{c_{d,q}}{c_d}$ converges to a normal distribution for $d
\rightarrow \infty$. Lastly, we consider in Section \ref{alternativeprobssection}
some related enumeration problems one might be interested in solving.

\acks
 Thanks to Adam Lawrence Epstein and Tan Lei for the suggestion of the problems and
to Christian Henriksen for a simplification of the method for finding the
generating functions. The author would furthermore like to thank Bodil Branner, Carsten Lunde Petersen, Carsten Thomassen, and again those mentioned
above for helpful discussions and
comments. The author is finally indebted to Philippe Flajolet for clarifying
comments about his and his co-authors' work on analytic combinatorics. \par
This research was supported by the European Union research training network CODY
(Conformal Structures and Dynamics).
\section{Background and Definitions}
 We present now a summary of some necessary concepts and definitions, and we include a new description of the combinatorial data set; it is different from the one presented in \cite{BD09} (it is in fact more similar to the definition of the combinatorial invariant in \cite{Sent}), but it is much more useful for our purposes here. For further details, please see \cite{BD09} and \cite{Sent}.\par
\begin{figure}%
    \centering
    \resizebox{!}{5cm}{\input{trajsatinfty.pstex_t}}
 \caption{The point at $\infty$ is a pole of order $d-2$ for vector fields $\xi_P \in \Xi_d$. There are $2(d-1)$ trajectories $\gamma_{\ell}$ which meet at infinity with asymptotic angles $  \frac{2 \pi\ell}{2\left(d-1\right)}$, $\ell \in \{ 0,1,\dots,2d-3 \}$. When the labelling index $\ell$ is even,  the trajectories are called \emph{incoming} to $\infty$, and when the index $\ell$ odd, they are called \emph{outgoing} from $\infty$.}
 \label{trajsatinfty}
 \end{figure}
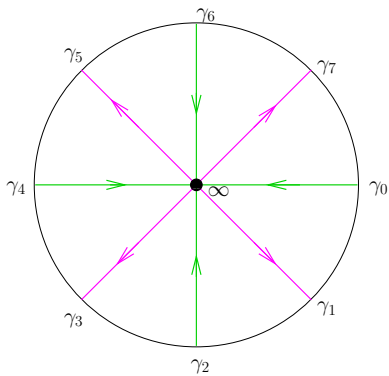

It can be shown that $\infty$ is a pole of order $d-2$ for vector fields $\xi_P \in \Xi_d$. There are $2(d-1)$ trajectories $\gamma_{\ell}$ which meet at infinity with asymptotic angles $  \frac{2 \pi\ell}{2\left(d-1\right)}$, $\ell \in \{ 0,1,\dots,2d-3 \}$. When the labelling index $\ell$ is even,  the trajectories are called \emph{incoming} to $\infty$, and when the index $\ell$ odd, they are called \emph{outgoing} from $\infty$ (see Figure \ref{trajsatinfty}).\par
\begin{figure}%
    \centering
    \resizebox{!}{5cm}{\input{trajsatinfty2.pstex_t}}
    \caption{There are $2d-2$ accesses to $\infty$ defined by the trajectories at infinity. An \emph{end} $e_{\ell}$ is infinity with access between $\gamma_{\ell-1}$ and $\gamma_{\ell}$. An \emph{odd end} is an end $e_k$ labelled by an odd index $k$, and an \emph{even end} is an end $e_j$ labelled by an even index $j$.}
    \label{trajsatinfty2}
 \end{figure}
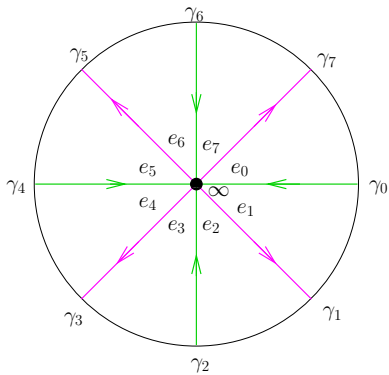

There are $2d-2$ accesses to $\infty$ defined by the trajectories at infinity. An \emph{end} $e_{\ell}$ is infinity with access between $\gamma_{\ell-1}$ and $\gamma_{\ell}$ (see Figure \ref{trajsatinfty2}). An \emph{odd end} is an end $e_k$ labelled by an odd index $k$, and an \emph{even end} is an end $e_j$ labelled by an even index $j$.\par
\emph{Separatrices}  $s_{\ell}$  are the maximal trajectories of $\xi_P$ incoming to and outgoing from $\infty$. They are labelled by the $2(d-1)$ asymptotic angles $  \frac{2 \pi\ell}{2\left(d-1\right)}$, $\ell \in \{ 0,1,\dots,2d-3 \}$.
  \begin{figure}%
    \centering
    \resizebox{!}{5cm}{\input{separatrixgraph2.pstex_t}}
    \caption{\emph{Separatrices}  $s_{\ell}$  are the maximal trajectories of $\xi_P$ incoming to and outgoing from $\infty$. They are labelled by the $2(d-1)$ asymptotic angles $  \frac{2 \pi\ell}{2\left(d-1\right)}$, $\ell \in \{ 0,1,\dots,2d-3 \}$. The separatrix structure can be
equivalently represented in a
\emph{separatrix disk model}, by labelling the points $\exp(\frac{2\pi {\rm
i} \ell}{2d-2})$, $\ell=0,\dots,2d-3$ on $\mathbb{S}^1$ by $s_{\ell}$ and
joining the points in the same equivalence class. In the figure, there are two homoclinic separatrices $s_{3,0}$ and $s_{7,6}$, and four landing separatrices $s_1$, $s_2$, $s_4$, and $s_5$. A few extra trajectories (in black) have been drawn in to better see the structure of the vector field. In this disk model, there are two sepal cells, one $\alpha \omega$ cell, and one center cell.}
    \label{separatrixgraph2}
 \end{figure}
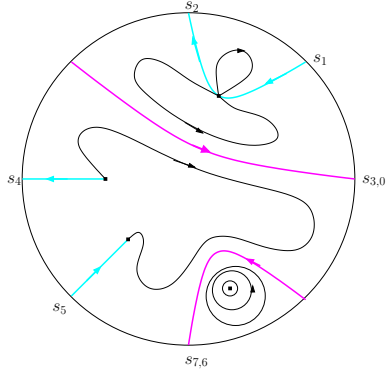
A separatrix $s_{\ell}$ is called \emph{landing} if $\bar{s}_{\ell}\setminus s_{\ell}=\zeta$, where $\zeta$ is an equilibrium point for $\xi_P$. A separatrix
$s_{\ell}=s_{k,j}$ is called \emph{homoclinic} if $\bar{s}_{k,j}\setminus s_{k,j}=\emptyset$ (see Figure \ref{separatrixgraph2}). A separatrix for a vector field $\xi_P \in \Xi_d$ can only be either homoclinic or landing. A homoclinic separatrix $s_{k,j}$ is labelled by the one odd index $k$ and the one even index $j$ corresponding to its two asymptotic directions at infinity.
The \emph{Separatrix graph:} $\Gamma_P=\bigcup \limits_{\ell=0}^{2d-3}\hat{s}_{\ell}$ completely determines the topological structure of the trajectories of a vector field (see, for instance, \cite{DN1975}, \cite{ALGM1973}).\par
 In \cite{BD09}, the separatrix structure is encoded via an equivalence relation and a marked subset $H$ on $\mathbb{Z}/(2d-2)$, where, in short, $H$ tells which separatrices are homoclinic, and for the remaining separatrices (which are landing), the equivalence relation essentially tells which $s_{\ell}$ land at the same equilibrium point (the formal definition will be given in Subsection \ref{combdatasetsection}). \par
 This separatrix structure can be
equivalently represented in a
\emph{separatrix disk model}, by labelling the points $\exp(\frac{2\pi {\rm
i} \ell}{2d-2})$, $\ell=0,\dots,2d-3$ on $\mathbb{S}^1$ by $s_{\ell}$ and
joining the points in the same equivalence class (see Figure \ref{separatrixgraph2}). \par
Conversely, if we define an abstract \emph{combinatorial data set}, $(\sim,H)$, which consists of a non-crossing
equivalence relation  $\sim$ and a marked subset $H$ on $\mathbb{Z}/(2d-2)$  which satisfy certain
properties (to be discussed), then it is known that there exists a monic and centered complex polynomial vector field of degree $d$ whose separatrix structure matches the data set \cite{BD09}.
\subsection{Zones}
The connected components $Z$ of $\mathbb{C}\setminus \Gamma_P$ are called \emph{zones}.
There are three types of zones for vector fields in $\Xi_d$, and the types of zones are determined by the types of their boundaries:
 \begin{itemize}
 \item
 A \emph{center zone} $Z$ contains an equilibrium point, which is a center, in its interior. Its boundary consists of one or several homoclinic separatrices and the point at infinity. If a center zone is on the left of $n$ homoclinic separatrices $s_{k_1,j_1},\dots,s_{k_n,j_n}$ on the boundary $\partial Z$, then the center zone has $n$ odd ends $e_{k_1},\dots,e_{k_n}$ at infinity on $\partial Z$ and the zone is called either a \emph{counter-clockwise center zone} or an \emph{odd center zone}. If a center zone is on the right of $n$ homoclinic separatrices $s_{k_1,j_1},\dots,s_{k_n,j_n}$ on the boundary $\partial Z$, then the center zone has $n$ even ends $e_{j_1},\dots,e_{j_n}$ at infinity on $\partial Z$ and the zone is called either a \emph{clockwise center zone} or an \emph{even center zone} (see Figure \ref{centerzoneshade1}).
  \begin{figure}%
    \centering
    \resizebox{!}{5cm}{\input{centerzoneshade1.pstex_t}}
    \caption{Pictured are the trajectories of a vector field with four center zones: one odd center zone (shaded) with homoclinic separatrices $s_{5,0}$, $s_{1,2}$, and $s_{3,4}$ and ends $e_1$, $e_3$, and $e_5$ on the boundary; and three even center zones, each with one homoclinic separatrix and one end on the boundary. }
    \label{centerzoneshade1}
 \end{figure}
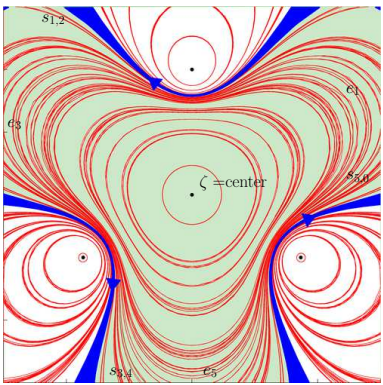
 \item
 A \emph{sepal zone} $Z$ has exactly one equilibrium point on the boundary, which is both the $\alpha$-limit point and $\omega$-limit point for all trajectories in $Z$ (i.e. $\zeta_{\alpha}= \zeta_{\omega}$). This equilibrium point is necessarily a multiple equilibrium point. The boundary $\partial Z$ contains exactly one incoming and one outgoing landing separatrix, and possibly one or several homoclinic separatrices, and the point at infinity.   If a sepal zone is to the left of $n$ homoclinic separatrices
$s_{k_1,j_1},\dots,s_{k_n,j_n}$ on its boundary, then it has $n+1$ odd ends on the
boundary: $e_{k_1},\dots,e_{k_n}$ and $e_{j_i+1}$ for some corresponding $j_i$,
depending on how one orders the separatrices. In this case, it is called an \emph{odd} sepal zone. Similarly, if a sepal zone is on the right of
$n$ homoclinic separatrices $s_{k_1,j_1},\dots,s_{k_n,j_n}$ on
its boundary, then it has $n+1$ even ends on the boundary, $e_{j_1},\dots,e_{j_n}$
and $e_{k_i+1}$ for some corresponding $k_i$, again depending on the ordering of
the separatrices. In this case, it is called an \emph{even} sepal zone (see Figure \ref{sepalzoneshade2}).
  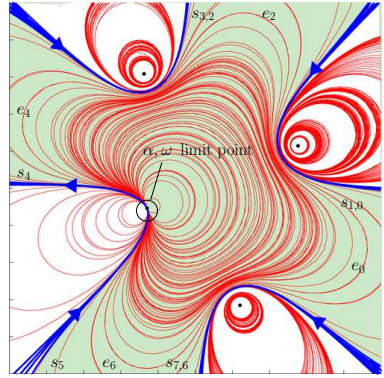
\begin{figure}%
    \centering
    \resizebox{!}{5cm}{\input{sepalzoneshade2.pstex_t}}
    \caption{Pictured are the trajectories of a vector field with an even sepal zone (shaded). On the boundary of the sepal zone is the double equilibrium point which is both the $\alpha$ and $\omega$ limit point of the trajectories; one incoming landing separatrix $s_4$ and one outgoing landing separatrix $s_5$; three homoclinic separatrices $s_{1,0}$, $s_{3,2}$, and $s_{7,6}$; and four ends at infinity $e_0$, $e_2$, $e_4$, and $e_6$. There is an odd sepal zone (not shaded) which shares the equilibrium point and the landing separatrices with the shaded sepal zone, but it has no homoclinic separatrices and only one odd end $e_5$ on the boundary.}
    \label{sepalzoneshade2}
 \end{figure}
\item
An \emph{$\alpha \omega$-zone} $Z$ has two equilibrium points on the boundary, $\zeta_{\alpha}\neq \zeta_{\omega}$, the $\alpha$-limit point and $\omega$-limit point for all trajectories in $Z$. The boundary $\partial Z$ contains one or two incoming landing separatrices and one or two outgoing landing separatrices, possibly one or several homoclinic separatrices, and the point at infinity. If an ${\rm \alpha \omega}$-zone is both on the left of $n_1$ homoclinic separatrices
$s_{k_1,j_1},\dots,s_{k_{n_1},j_{n_1}}$ and on the right of $n_2$ homoclinic
separatrices $s_{k_1,j_1},\dots,s_{k_{n_2},j_{n_2}}$ on the boundary, then the ${\rm
\alpha \omega}$-zone has $n_1+1$ odd ends ($e_{k_1},\dots,e_{k_{n_1}}$ and
$e_{j_i+1}$ for some corresponding $j_i$) and $n_2+1$ even ends
($e_{j_1},\dots,e_{j_{n_2}}$ and $e_{k_i+1}$ for some corresponding $k_i$) on the
boundary (see Figure \ref{alphaomegashade2}).\par
  \begin{figure}%
    \centering
    \resizebox{!}{5cm}{\input{alphaomegashade2.pstex_t}}
    \caption{Pictured are the trajectories of a vector field with an $\alpha \omega$-zone (shaded). On the boundary of the  zone are two equilibrium points: one which is the $\alpha$-limit point of the trajectories and the other is the $\omega$-limit point of the trajectories. Also on the boundary are one incoming landing separatrix $s_0$ and one outgoing landing separatrix $s_3$. The zone  is to the left of the homoclinic separatrix $s_{1,2}$  and on the right of the two  homoclinic separatrices $s_{5,4}$, and $s_{7,6}$. Finally, the boundary contains two odd ends $e_1$ and $e_3$ and three even ends $e_0$, $e_4$, and $e_6$ at infinity.}
    \label{alphaomegashade2}
 \end{figure}
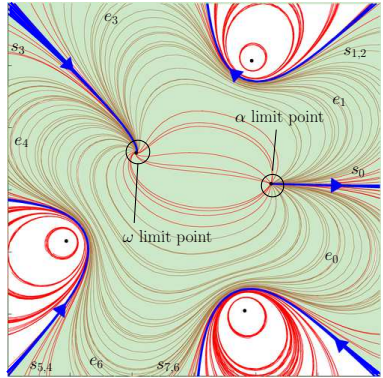
 \end{itemize}
 \proc{Remark.}
 It will be important to note that for a sepal zone, there is exactly one end whose index is not equal to an index of a homoclinic separatrix (in the notation above, $e_{j_i+1}$ for an odd sepal zone and $e_{k_i+1}$ for an even sepal zone).  Similarly, for an $\alpha \omega$-zone, there are exactly one odd end and one even end, neither of whose indices coincide with any index of a homoclinic separatrix (in the notation above the odd and even ends are $e_{j_i+1}$ and $e_{k_i+1}$ respectively).
 \medbreak
 In the separatrix disk model, the connected components of $\mathbb{D}\setminus \bigcup \limits_{\ell}s_{\ell}$ are called \emph{cells}. Cells are the abstract analog of zones, and are named accordingly (i.e. there are center cells, sepal cells, and $\alpha \omega$ cells). The types of cells are determined by the types of their boundaries (see Figure \ref{separatrixgraph2}).
 \subsection{Combinatorial Data Set}
 \label{combdatasetsection}
Now that we have defined the types of zones and corresponding cells, we can state the definition of the combinatorial data set as presented in \cite{BD09}.
\proc{Definition.}
\label{BDcombdatasetdef}
A \emph{combinatorial data set} $\left(\sim,H\right)$ of degree $d \geq 2$ consists of an equivalence relation $\sim$ on $\mathbb{Z}/\left(2d-2\right)$ and a marked subset $H \subset \mathbb{Z}/\left(2d-2\right)$ satisfying:
\begin{itemize}
\item[1)]
$\sim$ is non-crossing.
\item[2)]
If $\ell' \neq \ell''$, then $\ell' \sim \ell''$ and $\ell' \in H \Leftrightarrow \ell'' \in H$ and $\ell'$ and $\ell''$ have different parity.
\item[3)]
Every cell in the disk-model realization of $\left(\sim,H\right)$ is one of the five types: an $\alpha \omega$-cell, an odd or even sepal-cell, or an odd or even center-cell characterized as above.
\end{itemize}
\medbreak
\subsection{Transversals}
There is an equivalent way to encode the combinatorial structure of a vector field, which is much more useful for our purposes here. We define in this section the important structures needed to understand this equivalent definition of a combinatorial data set.\par
\begin{figure}%
    \centering
    \begin{minipage}{5cm}
    \resizebox{!}{5cm}{\input{transversals2_defensepres.pstex_t}}
    \end{minipage}
    \hspace{.25cm}
    \begin{minipage}{1cm}
    \resizebox{!}{1cm}{\input{rectmap.pstex_t}}
    \end{minipage}
    \hspace{.25cm}
    \begin{minipage}{2cm}
    \resizebox{!}{2cm}{\input{transversals_diststrip.pstex_t}}
    \end{minipage}
    \caption{There may be several transversal which avoid the equilibrium points and separatrices (the dashed curves), but there is exactly one distinguished transversal for each $\alpha \omega$-zone (in this case, $T_{3,0}$).}
    \label{transstrip}
 \end{figure}
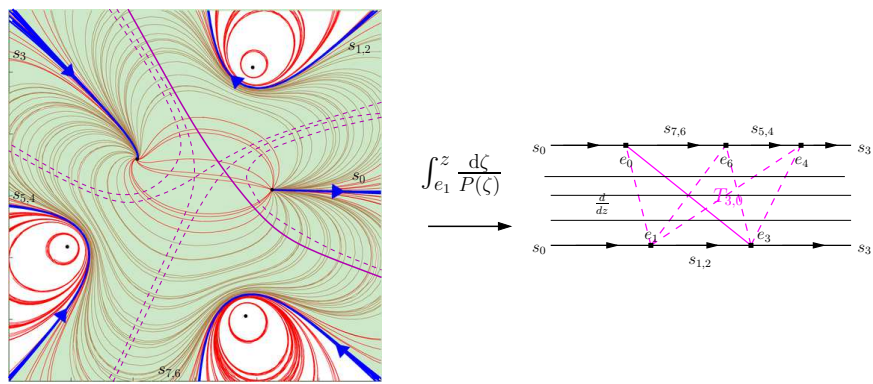
\begin{figure}
\centering
\mbox{\resizebox{!}{6cm}{\input{disttransvf.pstex_t}}\quad \resizebox{!}{3.5cm}{\input{disttrans.pstex_t}}}
\caption{Each $\alpha \omega$-zone is isomorphic to a strip. We define the distinguished transversal to be the geodesic in the metric $|dz|/|P(z)|$ joining the ends $e_k$ and $e_j$ (in this figure, $e_3$ and $e_0$) where $e_j$ is the left-most end on the upper boundary of the strip and $e_k$ is the right-most end on the lower boundary of the strip. Since the indices of the distinguished transversal are the same as the indices for the two landing separatrices on the upper left and lower right boundary of the strip, they can never coincide with the indices for a homoclinic separatrix.}
\label{disttrans}
\end{figure}
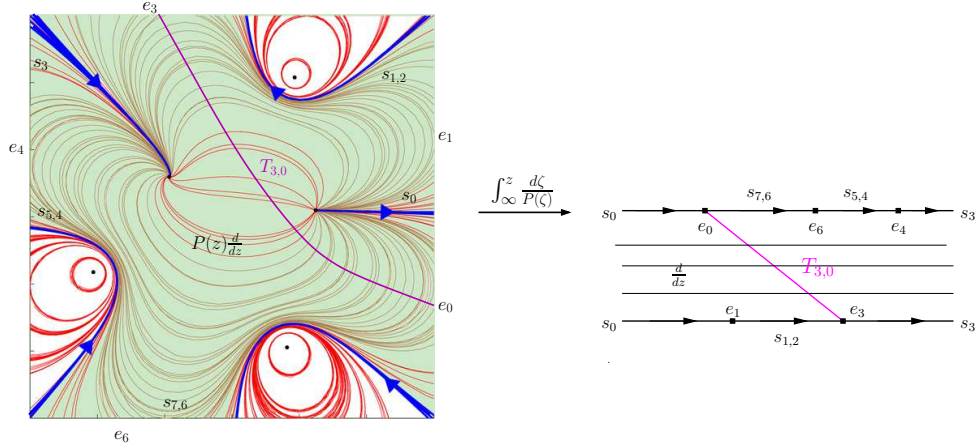

In any simply connected domain avoiding zeros of $P$, the differential $\frac{\rm{d}z}{P(z)}$ has an antiderivative,
unique up to addition by a constant
\begin{equation}
\phi(z)=\int_{z_0}^{z} \frac{\rm{d}w}{P(w)}.\nonumber
\end{equation}
Note that
\begin{equation}
\phi_{\ast}\left(\xi_P\right)=\phi'\left(z\right)P\left(z\right)\frac{{\rm d}}{{\rm d}z}=\frac{{\rm d}}{{\rm d}z}.
\end{equation}
The coordinates $w=\phi(z)$ are, for this reason, called \emph{rectifying coordinates}.
We will call the images of zones under rectifying coordinates \emph{rectified zones}. The rectified zones are of the following types:
\begin{itemize}
\item
The image of an $\alpha \omega$-zone under $\phi$ is a horizontal strip (see Figure \ref{alphaomegashade}).
  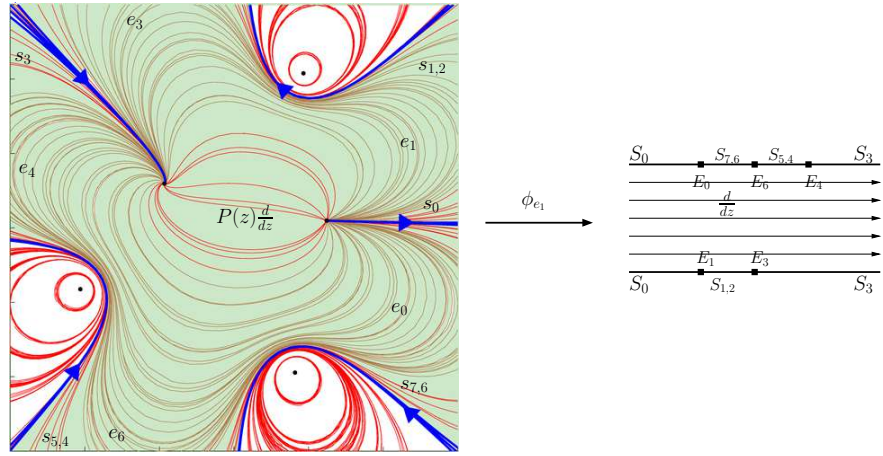
\begin{figure}%
    \centering
    \resizebox{!}{6cm}{\input{alphaomegashade.pstex_t}}
    \caption{The image of an $\alpha \omega$-zone under $\phi$ is a horizontal strip.}
    \label{alphaomegashade}
 \end{figure}
\item
The image of an odd sepal zone under $\phi$ is an upper half plane, and the image of an even sepal zone is a lower half plane (see Figure \ref{sepalzoneshade}).
  \begin{figure}%
    \centering
    \resizebox{!}{6cm}{\input{sepalzoneshade.pstex_t}}
    \caption{The image of an odd sepal zone under $\phi$ is an upper half plane, and the image of an even sepal zone is a lower half plane. In this figure, there is an even sepal zone mapped to a lower half plane.}
    \label{sepalzoneshade}
 \end{figure}
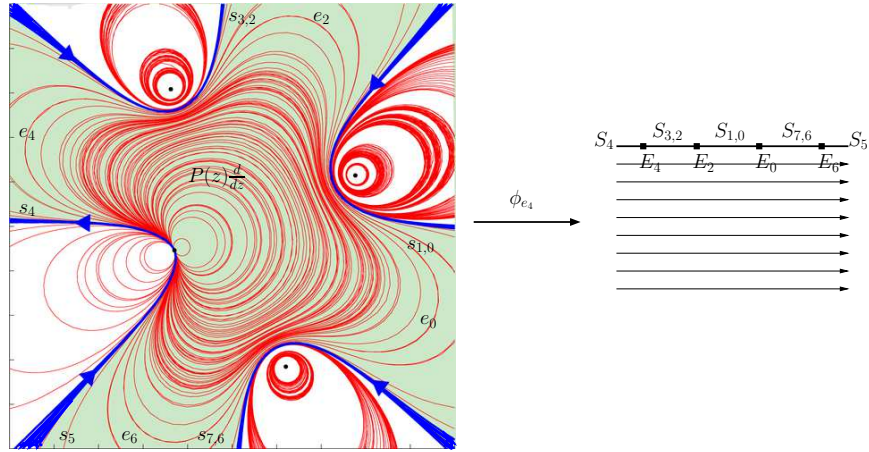
\item
 The image of a center zone (minus a curve contained in the zone which joins the center $\zeta$ and $\infty$) under $\phi$ is a vertical half strip. It is an upper vertical half strip for an odd center zone and a lower vertical half strip for an even center zone (see Figure \ref{centerzoneshade}).\par
  \begin{figure}%
    \centering
    \resizebox{!}{6cm}{\input{centerzoneshade.pstex_t}}
    \caption{The image of a center zone (minus a curve contained in the zone which joins the center $\zeta$ and $\infty$) under $\phi$ is a vertical half strip. It is an upper vertical half strip for an odd center zone and a lower vertical half strip for an even center zone. In this case, there is an odd center zone mapped to an upper vertical half strip.}
    \label{centerzoneshade}
 \end{figure}
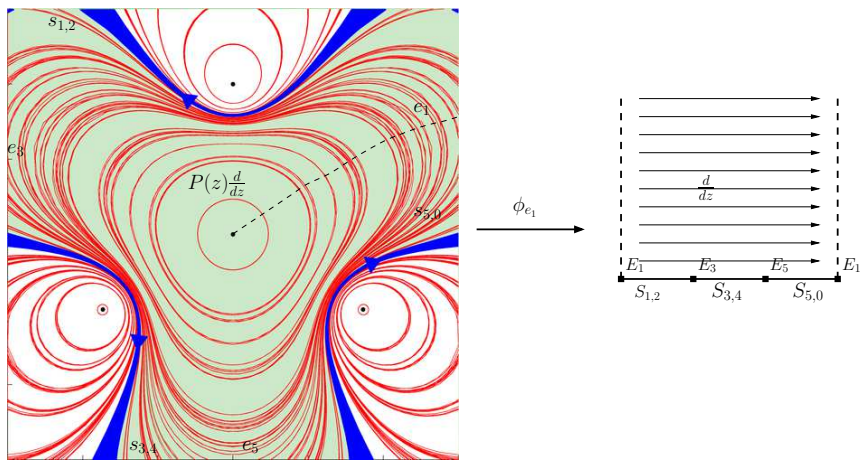
\end{itemize}
 Via the rectifying coordinates, it is evident that there are a number of closed geodesics in $\hat{\mathbb{C}} \setminus \{\text{equilibrium pts}  \}$ in the metric with length element $\frac{|\rm{d}z|}{|P(z)|}$ through $ \infty$. Among these are the $h$ homoclinic separatrices, and there are $s$ \emph{distinguished transversals} (defined below).
 \proc{Definition.}
  The \emph{distinguished transversal} $T_{k,j}$ is the geodesic in the metric $\frac{|\rm{d}z|}{|P(z)|}$ joining the ends $e_k$ and $e_j$, avoiding the separatrices and equilibrium points, where $e_j$ is the left-most end on the upper boundary and $e_k$ is the right-most end on the lower boundary of the strip which is the image of the $\alpha \omega$-zone the transversal is contained in (see Figure \ref{transstrip} and \ref{disttrans}).
 \medbreak
 Note that the way in which the distinguished transversal is chosen, the indices of the ends it joins are exactly those ends whose indices will never coincide with the indices of any homoclinic separatrices.
\subsection{Alternative Combinatorial Description}
We will prove that one can equivalently describe the combinatorics by the union of homoclinic separatrices $s_{k,j}$ and distinguished transversals $T_{k,j}$.
We state the definition below, and prove that it is equivalent to the one presented before (from \cite{BD09}).
 \proc{Definition.}
A \emph{combinatorial data set} $(\iota,H,T)$ of degree $d\geq 2$ consists of an involution $\iota : \mathbb{Z}/(2d-2)\rightarrow \mathbb{Z}/(2d-2)$ and marked subsets $H$ and $T$ satisfying:
\begin{itemize}
\item
If $\iota (\ell) \neq \ell$, then $\ell$ and $\iota(\ell)$ have different parity (if one is even, then the other is odd).
\item
The involution $\iota$ is non-crossing.
\item
If $\iota (\ell) \neq \ell$, then either both $\ell$ and $\iota(\ell)$ are in $H$ or both $\ell$ and $\iota(\ell)$ are in $T$.
\end{itemize}
 \medbreak
  This definition leads to a description of the combinatorics via non-crossing pairings of even and odd numbers.
 The pairing of indices in $H$ represents pairings of indices $k$ and $j$ corresponding to homoclinic separatrices $s_{k,j}$, and the pairing of indices in $T$ represents pairings of indices $k$ and $j$ corresponding to indices of ends $e_k$ and $e_j$ that the distinguished transversal $T_{k,j}$ joins. The indices which are not paired (i.e. $\ell$ such that $\iota(\ell)=\ell$) correspond to the indices of ends $e_{\ell}$ on the boundary of sepal zones which do not coincide with indices of any homoclinic separatrices. Essentially, we want to use the numbers $\mathbb{Z}/(2d-2)$ to stand for indices of separatrices for homoclinics, and indices of ends otherwise.  Since these were chosen in the way described above, this never causes a conflict.  \par
 As one can draw a disk model with a separatrix graph given the first definition of combinatorial invariant, one can draw a transversal graph given the second definition. The transversal structure can be
equivalently represented in a
\emph{transversal disk model}, by labelling the points $\exp(\frac{2\pi {\rm
i} \ell}{2d-2})$, $\ell=0,\dots,2d-3$ on $\mathbb{S}^1$ by $s_{\ell}$ and the points between them $\exp(\frac{2\pi {\rm
i} (\ell-1/2)}{2d-2})$ by $e_{\ell}$, joining $s_{k}$ and $s_j$ corresponding to homoclinic separatrices $s_{k,j}$, and
joining $e_k$ and $e_j$ corresponding to each distinguished transversal $T_{k,j}$ (see Figure \ref{transversalgraph}). That the involution $\iota$ is non-crossing means that the transversal and homoclinic separatrices in the transversal disk model do not cross. \par
  \begin{figure}%
    \centering
    \resizebox{!}{6cm}{\input{transversalgraph.pstex_t}}
    \caption{Pictured is the transversal disk model equivalent to the separatrix disk model pictured in Figure \ref{separatrixgraph2}. There are two homoclinic separatrices $s_{3,0}$ and $s_{7,6}$, and there is one distinguished transversal $T_{5,4}$. Note that each connected component contains exactly one equilibrium point. }
    \label{transversalgraph}
 \end{figure}
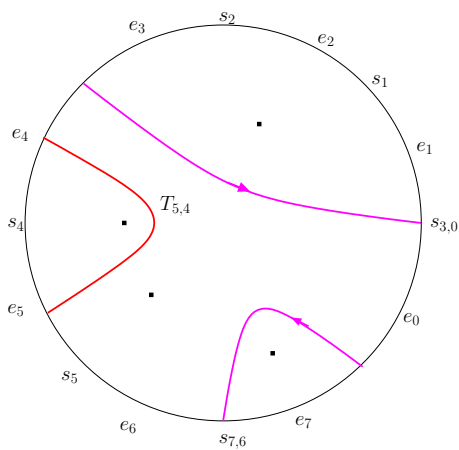
 \begin{prop}
 The definition above is equivalent to the definition of combinatorial data set presented in \cite{BD09} (Definition \ref{BDcombdatasetdef})
 \end{prop}
 \proc{Proof.}
 We basically need to prove that separatrix disk models are in one-to-one correspondence with transversal disk models. \par
 If we start with a separatrix disk model, such that every cell is one of the five types, then we can construct a unique transversal disk model.  Note first that the homoclinic separatrices stay the same. Then for each $\alpha \omega$-zone, there are exactly one even end $e_j$ and one odd end $e_k$ whose indices do not coincide with any of the homoclinic indices on the boundary. Join these ends in the disk. Finally, remove all landing separatrices from the separatrix disk model, and what you are left with is a transversal disk model (see Figure \ref{transversalgraph}). \par
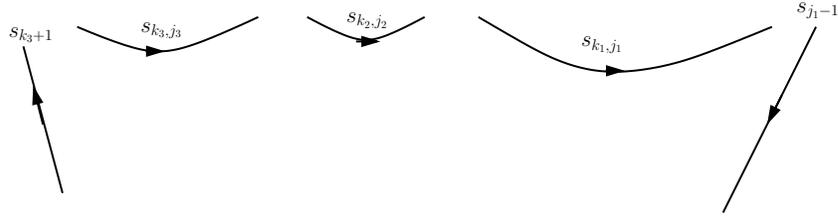
\begin{figure}%
    \centering
    \resizebox{!}{3cm}{\input{openHchains1.pstex_t}}
    \caption{An open chain of three homoclinic separatrices $s_{k_1,j_1},\ s_{k_2,j_2}$, and $s_{k_{3},j_{3}}$ such that $j_{i+1}=k_{i}+1$ for $i=1,2$ and $s_{j_1-1}$ and $s_{k_3+1}$ are landing. }
    \label{openHchains1}
 \end{figure}

 Now conversely given a valid transversal disk model (i.e. one satisfying the definition above), we again leave the homoclinic separatrices alone as before.  The homoclinic separatrices and distinguished transversals decompose the transversal disk model into $h+s+1$ connected components. Call each of these components a transversal cell (as not to confuse them with connected components with respect to a separatrix disk model). Let all $s_{\ell}$ on the boundary of each transversal cell be joined to the same point in the interior of the transversal cell. If there are no $s_{\ell}$ on the boundary of a transversal cell, then that component is a center cell (with respect to the separatrix disk model).  Now remove the curves corresponding to the distinguished transversals. This procedure gives a separatrix disk model, and now we just need to show that this separatrix disk model has the allowable types of cells. As already mentioned, if there are only homoclinic separatrices on the boundary of a separatrix cell, then that component is a center cell.  If there is an open chain of $n$ (where we allow $n=0$) homoclinic separatrices: $s_{k_1,j_1},\dots,s_{k_{n},j_{n}}$ such that either $j_{i+1}=k_{i}+1$ for $i=1,\dots,n-1$ and $s_{j_1-1}$ and $s_{k_n+1}$ are landing (see Figure \ref{openHchains1}) or if $k_{i+1}=j_{i}+1$ for $i=1,\dots,n-1$ and $s_{k_1-1}$ and $s_{j_n+1}$ are landing (see Figure \ref{openHchains2}) on the boundary of the same separatrix cell, then either the two landing separatrices land at the same point, or they do not.  If they do land at the same point, then $s_{j_1-1}$ and $s_{k_n+1}$  (or $s_{k_1-1}$ and $s_{j_n+1}$) are the two landing separatrices on the boundary of a sepal cell. If they land at different points, then they must be separated by a distinguished transversal $T_{k,k_n+1}$ for some $k$ associated to another open homoclinic chain (or $T_{j_n+1,j}$ for some $j$ associated to another open homoclinic chain) on the boundary of the same separatrix cell. Then $s_{j_1-1}$ must land at the same point as $s_k$, since they are on the boundary of the same connected component in the transversal disk model (similarly for $s_{k_1-1}$ and $s_j$) (see Figure \ref{openHchains}). One can make a similar argument for the other open homoclinic chain involved, and then one must get an $\alpha \omega$-cell.
 \ep
 \medbreak
  \begin{figure}%
    \centering
    \resizebox{!}{3cm}{\input{openHchains2.pstex_t}}
    \caption{An open chain of three homoclinic separatrices $s_{k_1,j_1},\ s_{k_2,j_2}$, and $s_{k_{3},j_{3}}$ such that $k_{i+1}=j_{i}+1$ for $i=1,2$ and $s_{k_1-1}$ and $s_{j_n+1}$ are landing. }
    \label{openHchains2}
 \end{figure}
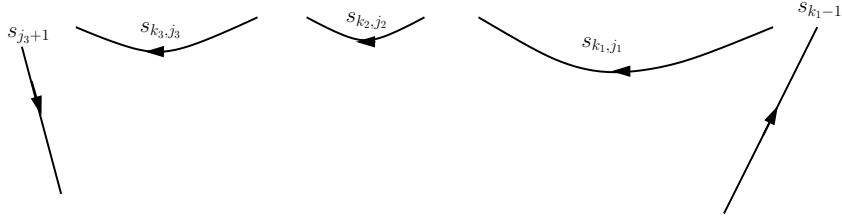
  \begin{figure}%
    \centering
    \resizebox{!}{4.5cm}{\input{openHchains.pstex_t}}
    \caption{An open chain of three homoclinic separatrices $s_{k_1,j_1},\ s_{k_2,j_2}$, and $s_{k_{3},j_{3}}$ such that $k_{i+1}=j_{i}+1$ for $i=1,2$ and $s_{k_1-1}$ and $s_{j_n+1}$ are landing. If they land at different points, then they must be separated by a distinguished transversal $T_{j_3+1,j}$ for some $j$ associated to another open homoclinic chain on the boundary of the same separatrix cell. Then $s_{k_1-1}$  must land at the same point as $s_j$, since they are on the boundary of the same connected component in the transversal disk model.}
    \label{openHchains}
 \end{figure}
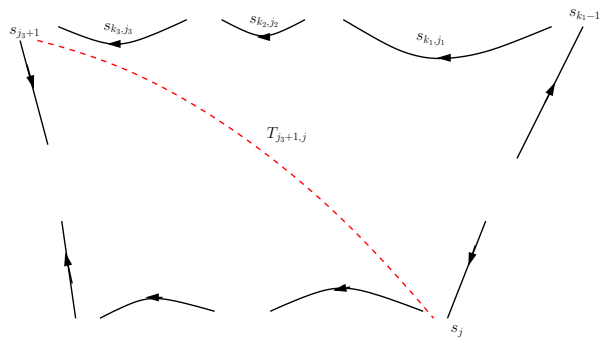
\section{Setup of the Bracketing Problem}
\label{bracketingsetupsection}
 We have seen that the combinatorial structure of a monic and centered complex polynomial vector field can be completely described by pairings of even integers with odd integers, and  conversely that any such pairings as defined by the transversal combinatorial data set leads to the combinatorics for a monic and centered complex polynomial vector field.
The goal is to utilize the alternative combinatorial description in order to convert the problem of
counting the
combinatorial classes into a so-called \emph{bracketing problem}: a
combinatorial problem involving pairings of parentheses placed in a string of
elements (see \cite{Com}) in a \emph{valid} way (to be defined). \par
In short, each homoclinic separatrix $s_{k,j}$ has exactly one even and one odd index associated with it, and each distinguished transversal $T_{k,j}$ has exactly one even and one odd index associated to it. Therefore, each combinatorial data set can be fully described by pairings of even and odd numbers, corresponding to these mutually disjoint pairs of indices, and hence they can be counted by a bracketing problem. In our case, we will have balanced parenthetical configurations placed in the string $0 \ 1 \ 2\dots 2d-3$, representing the indices in $\mathbb{Z}/(2d-2)$, and we will use round brackets $(\dots)$ to pair an even an odd number corresponding to homoclinics and square brackets $[\dots]$ to pair an even and odd number corresponding to distinguished transversals. We now elaborate on these ideas.
\subsection{Structurally Stable Vector Fields}
\label{strstablecase}
For peadagogical reasons, we start by describing the case of counting
structurally stable complex polynomial vector fields by a bracketing problem.\par
\begin{theo}[Douady, Estrada, Sentenac]
The number of combinatorial classes for the structurally stable vector fields in $\Xi_d$ is the Catalan number
\begin{equation}
C_{d-1}=\frac{1}{d}\binom{2(d-1)}{d-1}.
\end{equation}
\end{theo}
\proc{Proof.}
For structurally stable vector fields, the combinatorics is completely describable
by the $d-1$ transversals (since there are neither homoclinic separatrices nor
sepal zones), and hence non-crossing pairings of each odd end $e_k$ with an
even end $e_j$ (this is, in fact, how \cite{Sent} defines the combinatorial invariant for the structurally stable vector fields). This pairing of ends can be represented by placing parentheses
in the string $0 \ 1 \ 2\dots 2d-3$, where the elements paired by parentheses
correspond to the labels of the ends we want to pair.\footnote{A more familiar
description might
be evident in the disk model, where this pairing can be seen as the number of
non-crossing handshakes of $2(d-1)$ people seated around a round table.}
Note
that this is equivalent to the number of ways to make \emph{valid} pairings of
$d-1$ pairs of parentheses, where the string $0 \ 1 \ 2\dots 2d-3$ is not actually
necessary, since every end is paired with another. For example, if $d-1=3$, then the configurations in the table below are equivalent:
\begin{center}
\begin{tabular}{| c | c |}
\hline
[01][2[34]5]& [ ][ [ ] ] \\ \hline
[0[12]3][45] & [ [ ] ][ ] \\ \hline
[0[1[23]4]5] & [ [ [ ] ] ] \\ \hline
[0[12][34]5] & [ [ ][ ] ] \\ \hline
[01][23][45] & [ ][ ][ ] \\ \hline
\end{tabular}
\end{center}
You might recognize then that
the number of combinatorial classes (with the labelling of the separatrices) for
the structurally stable monic and centered complex polynomial vector fields of degree $d$ is just the Catalan number $C_{d-1}$. The proof is classical (see for instance \cite{Davis}), but we include it here for completeness.\par
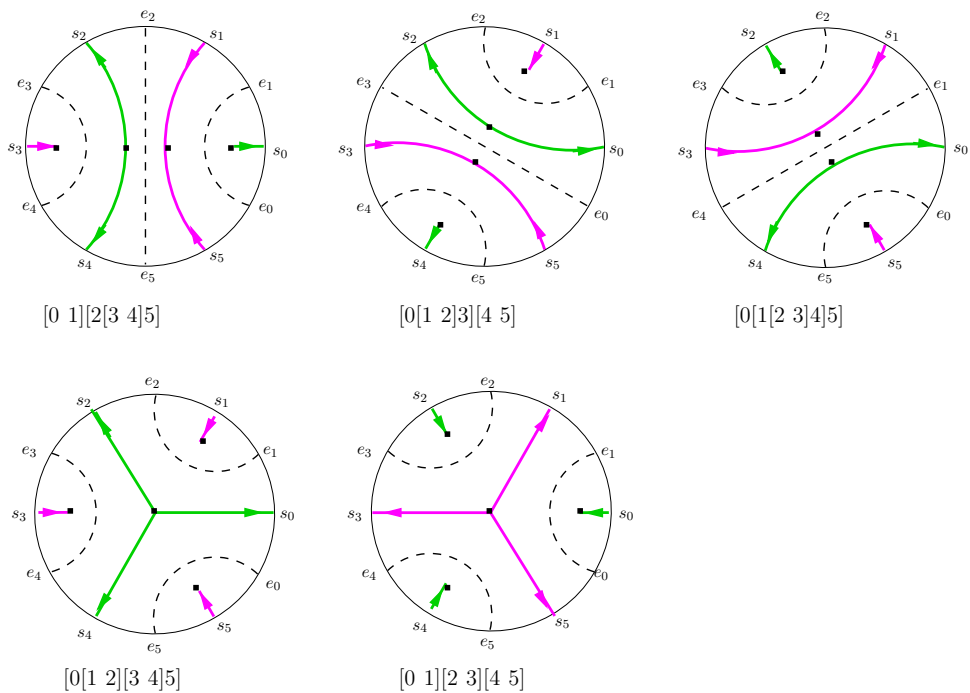
\begin{figure}
\resizebox{!}{9.2cm}{\input{d4generic.pstex_t}}
\caption{The five different disk models for the structurally stable vector fields of degree $d=4$.
The involution on the ends is marked by the dashed curves. The representation of
the involution in brackets is displayed below each figure.}
\label{d4generic}
\end{figure}

The goal is to write a recursion equation for $C_n$, the $n^{\rm th}$ Catalan number, which counts the number of configurations matching $n$ pairs of parentheses. The first character of any balanced configuration is an open parenthesis "[". Somewhere in the configuration is the matching "]" for the open one. In between that pair of parentheses is a balanced configuration, and to the right is another balanced configuration:
\begin{equation}
[A]B,
\end{equation}
where $A$ and $B$ are balanced parenthetical configurations. Each $A$ and $B$ can have anywhere from 0 to $n-1$ pairs of parentheses, but together they must have exactly $n-1$ pairs of parentheses. So if $A$ has $k$ pairs of parentheses, then $B$ must have $n-k-1$ pairs of parentheses. Thus we count all configurations where $A$ has 0 pairs and $B$ has $n-1$ pairs, $A$ has 1 pair and $B$ has $n-2$ pairs, etc. Add them up, and you get the total number of configurations with $n$ balanced pairs of parentheses. Therefore,
\begin{equation}
\label{catrecursion}
C_n=C_0C_{n-1}+C_1C_{n-2}+ \cdots +C_{n-2}C_{1}+C_{n-1}C_{0}, \quad C_0=1.
\end{equation}
We will use this recursion equation and the generating function
\begin{equation}
\label{catgenfunc}
f(z)=\sum \limits_{n=0}^{\infty}C_nz^n
\end{equation}
in order to derive a closed-form expression for $C_n$.  Using \eqref{catgenfunc}, one can compute
\begin{equation}
\left[f(z) \right]^{2}=C_0C_0+\left(C_0C_1+C_1C_0 \right)z+\left(C_0C_2+C_1C_1+C_2C_0 \right)z^2+\dots.
\end{equation}
Combining this with \eqref{catrecursion}, one arrives at
\begin{equation}
\left[f(z) \right]^{2}=C_1+C_2z+C_3z^2+C_4z^3+\dots.
\end{equation}
Therefore,
\begin{equation}
f(z)=C_0+z\left[f(z) \right]^{2}.
\end{equation}
Applying the quadratic formula,
\begin{equation}
f(z)=\frac{ 1 \pm \sqrt{1-4z} }{2z},
\end{equation}
and one must choose the branch
\begin{equation}
f(z)=\frac{ 1 - \sqrt{1-4z} }{2z},
\end{equation}
since $f(z) \rightarrow 1$ as $z \rightarrow 0$, which is desired since $f(z)=C_0=1$ (choosing the "+" branch gives $f(z) \rightarrow \infty$ as $z \rightarrow 0$). Expanding $\sqrt{1-4z}=\left(1-4z \right)^{1/2}$ using the binomial formula and making some arithmetic manipulations, one arrives at
\begin{equation}
f(z)=\sum \limits_{n=0}^{\infty} \frac{1}{n+1}\binom{2n}{n}z^n,
\end{equation}
giving the desired result.
\ep
\medbreak
An example is given in Figure \ref{d4generic} to demonstrate.
\subsection{Non-Structurally Stable Vector Fields}
 The convention, as stated before, is that round
parentheses $(\cdots)$ are used to mark pairings corresponding to homoclinic
separatrices, square parentheses $[\cdots]$ are used to mark pairings corresponding to a
distinguished transversal in each
$\alpha \omega$-zone, and elements that are not paired correspond to the
sepal-zones.  \par
\begin{figure}
\resizebox{!}{4.5cm}{\input{d4wsepal.pstex_t}}
\caption{Disk models for two examples of vector fields of degree $d=4$ having
sepal zones and no homoclinic separatrices.  The pairing of the ends is marked by
the dashed curves. The
representation of the combinatorics in brackets is displayed below each figure.}
\label{d4wsepal}
\end{figure}
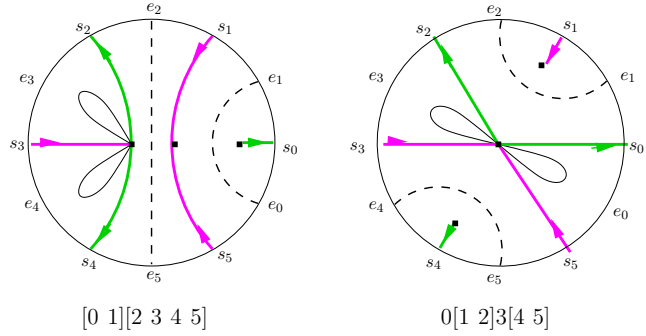

We first look at the case when there are sepal zones but no homoclinic separatrices. Therefore, the zones can only be either
$\alpha \omega$-zones, having exactly one even and one odd end, or sepal zones,
having exactly one even or one odd end. The string $0 \ 1 \ 2\dots 2d-3$ represents
the labels of the $2d-2$ ends in this case. The ${\rm\alpha \omega}$-zones will then
be denoted by square parentheses $[\cdots]$ in this string, and the labels corresponding to
the ends of sepal zones are not paired by a set of parentheses. We see now that the
string $0 \ 1 \ 2\dots 2d-3$ is  necessary since the unpaired elements will be
placeholders corresponding to the sepal zones. Two
examples are given in Figure \ref{d4wsepal}.\par
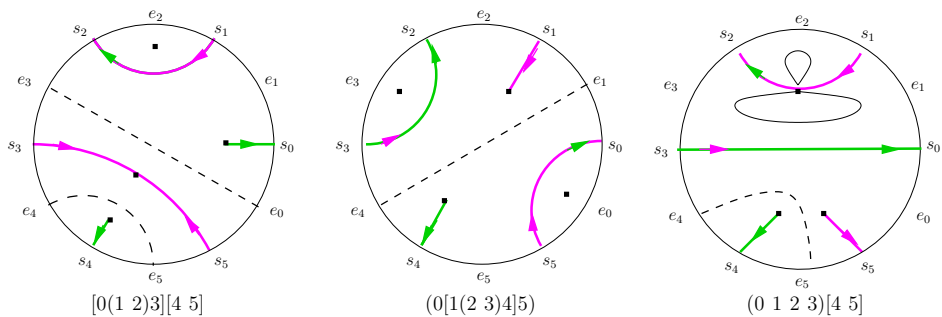
\begin{figure}
\resizebox{!}{4.2cm}{\input{d4wsepalnhom.pstex_t}}
\caption{Disk models for three examples of vector fields of degree $d=4$ having
sepal zones or/and homoclinic separatrices.  The pairing of the ends is marked by
the dashed curves. The representation of the combinatorics in brackets is
displayed below each figure.}
\label{d4wsepalnhom}
\end{figure}

What then happens if we introduce homoclinic separatrices? Then it is no longer
enough to consider pairings of ends. The string $0 \ 1 \ 2\dots 2d-3$ now represents
the indices which may belong to either the ends $e_{\ell}$ or separatrices $s_{\ell}$ (as was in the involution definition of the combinatorial data set). For every homoclinic
separatrix $s_{k,j}$, we pair the numbers $k$ and $j$ by round parentheses $(\cdots)$.
For example, $0 \ 1 \ \dots (j \dots k)\dots$ or $0 \ 1 \ \dots (k \dots
j)\dots$. As usual, for every distinguished transversal $T_{k,j}$, we pair the numbers $k$ and $j$ by square parentheses $[\cdots]$ (see Figure
\ref{d4wsepalnhom}).\par
We can see from this that given any combinatorial data set, there is a unique
bracketing representation defined in this way.
\subsection{Valid Bracketings}
The bracketing must satisfy certain rules so that they are in accordance with
what can happen for a vector field (and hence a combinatorial data set). Pairs
of square and round parentheses placed in a string of elements $0 \ 1 \ 2 \dots
2d-3$ is called a \emph{valid bracketing} if:
\begin{enumerate}
\item there are an equal number of right and left parentheses for each type,
\item the number of left parentheses must be greater than or equal to the number
of right, reading from left to right {\bf Example: ())()( not valid},
\item there must be at least one element between successive left (resp. right)
parentheses {\bf Example: $((01)23)$ not valid},
\item there must be an even number of elements in each pair of parentheses {\bf
Example: $(012)3$ not valid},
\item and the square brackets must not cross the round brackets {\bf Example:
([)] not valid}
\end{enumerate}
It is clear that a valid bracketing as defined above leads to an admissible combinatorial data set.
\subsection{Recursion Equations and Generating Functions}
\label{genfctsection}
We compute the recursion equation for $p_n$, the number of ways to place any
number of round or square parentheses in a valid way in a string of $n$ elements,
by using the following algorithm to generate any possible valid parenthetical
configuration in a string of $n$ elements, which we will denote by $s_n$.
\begin{equation}
s_n=\begin{cases}
"s"s_{n-1}&\text{if string does not start with parenthesis } \\
"(s"s_{2a} "s)" s_b&\text{if string starts with round parentheses }\\
"[s"s_{2a} "s]" s_b&\text{if string starts with square parentheses }\end{cases}
\end{equation}
for non-negative integers $a$ and $b$ satisfying $2a+b+2=n$, and the $s$ representing numbers in $\{ 0, 1 , \dots, 2d-3\}$. A table of the construction of such strings from $n=1$ to $n=5$ is given below. \newline
\begin{tabular}{| l | l | l | l| l | l |}
\hline
$n$ & $"s"s_{n-1}$ & $a$ & b & $"(s"s_{2a} "s)" s_b$ & $"[s"s_{2a} "s]" s_b$ \\ \hline
1 & $s$ & & &&\\ \hline
2 & $ss$ & 0 & 0&  $(ss)$& $[ss]$ \\
3& $sss$ & 0 & 1 &$(ss)s$  & $[ss]s$\\
& $s(ss)$ & & & & \\
& $s[ss]$ &  & & & \\ \hline
4 & $ssss$ & 0& 2 & $(ss)ss$ & $[ss]ss$ \\
& $ss(ss)$ & && $(ss)(ss)$ & $[ss](ss)$ \\
& $ss[ss]$ & &&$(ss)[ss]$ & $[ss][ss]$ \\
& $s(ss)s$ &1 & 0 & $(ssss)$ & $[ssss]$ \\
& $s[ss]s$ &&& $(s(ss)s)$ & $[s(ss)s]$ \\
&&&&$(s[ss]s)$ & $[s[ss]s]$ \\ \hline
5 & $sssss$& 0&3& $(ss)sss$& $[ss]sss$\\
& $sss(ss)$ & && $(ss)s(ss)$ & $[ss]s(ss)$ \\
& $sss[ss]$ & && $(ss)s[ss]$& $[ss]s[ss]$ \\
& $ss(ss)s$ & &  & $(ss)(ss)s$ &  $[ss](ss)s$\\
& $ss[ss]s$ &&&   $(ss)[ss]s$ & $[ss][ss]s$ \\
& $s(ss)(ss)$& 1 & 1 & $(ssss)s$ & $[ssss]s$ \\
& $s[ss](ss)$ & & & $(s(ss)s)s$ & $[s(ss)s]s$ \\
& $s(ss)[ss]$&&& $(s[ss]s)s$ & $[s[ss]s]s$ \\
& s[ss][ss]&&&&\\
& s(ssss)&&&&\\
& s[ssss]&&&&\\
& s(s(ss)s)&&&&\\
& s[s(ss)s]&&&&\\
& s(s[ss]s)&&&&\\
& s[s[ss]s]&&&&\\
\hline
\end{tabular}
\vspace{.25cm}\newline
Then $p_n$ satisfies the recursion equation
\begin{equation}
p_n=\begin{cases}
p_{n-1}+2\left( \sum \limits_{2a+b+2=n}p_{2a} p_b \right), \quad &n\geq 1 \text{ and } a, b \geq 0\\
1, \quad &n=0\\
0, \quad &n<0.\end{cases}
\end{equation}
We are interested in the generating function
\begin{equation}
G(z)=\sum \limits_{d=1}^{\infty}c_d z^d,
\end{equation}
where $c_d=p_{2(d-1)}$, since for degree $d$ the string we work with is $0 \dots 2d-3$. Note that $G(0)=0$. For $n \geq 0$, let $q_n=p_{2n}$ and
$r_n=p_{2n-1}$ ($r_0=p_{-1}=0$). Let also
\begin{equation}
f(z)=\sum \limits_{n=0}^{\infty}q_nz^{n}
\end{equation}
and
\begin{equation}
g(z)=\sum \limits_{n=1}^{\infty}r_nz^{n}.
\end{equation}
Note that $zf=G$. One can deduce for $n\geq 0$
\begin{align}
\label{qn}
q_n=p_{2n}&=p_{2n-1}+2\sum\limits_{j=0}^{n-1}p_{2j}p_{2n-2j-2}\nonumber \\
&=r_n+2\sum\limits_{j=0}^{n-1}q_j q_{n-j-1}
\end{align}
and
\begin{align}
\label{rn}
r_n=p_{2n-1}&=p_{2n-2}+2\sum\limits_{j=0}^{n-1}p_{2j}p_{2n-2j-3}\nonumber \\
&=q_{n-1}+2\sum\limits_{j=0}^{n-1}q_j r_{n-j-1}.
\end{align}
From equations \eqref{qn} and \eqref{rn}, one can deduce
\begin{equation}
f=1+g+2zf^2 \quad \text{and} \quad g=zf+2zfg
\end{equation}
respectively, giving
\begin{equation}
\label{algeq1}
4z^2f^3-4zf^2+(z+1)f-1=0.
\end{equation}
Substituting $zf=G$, one gets that $G$ satisfies the algebraic equation
\begin{equation}
G^3-G^2+\frac{(z+1)}{4}G-\frac{z}{4}=0.
\end{equation}
For the complete problem, namely giving the number of combinatorial classes of each dimension, we are interested in the generating function
\begin{equation}
G(z,t)=\sum\limits_{d=1,q=0}^{\infty}c_{d,q}z^dt^q, \quad c_{d,q}=p_{2(d-1),q},
\end{equation}
where $q$ is the real dimension of the combinatorial class in parameter space. Note that
\begin{equation}
G(z,1)=\sum\limits_{d=1,q=0}^{\infty}c_{d,q}z^d=G(z).
\end{equation}
Folllowing the same method as before,
\begin{equation}
s_{n,q}=\begin{cases}
"s"s_{n-1,q}&\text{if string does not start with parenthesis } \\
"(s"s_{2a,q_1} "s)" s_{b,q_2}& q_1+q_2+1=q\\
"[s"s_{2a,q_1} "s]" s_{b,q_2}& q_1+q_2+2=q\end{cases}
\end{equation}
for $2a+b+2=n$.
The condidtions on $q_1$ and $q_2$ come from the fact that each homoclinic separatrix (and hence each round pair of parentheses) contributes real dimension 1, and each $\alpha \omega$-zone (and hence each square pair of parentheses) contributes complex dimension 1 (real dimension 2). Manipulating the coefficients as in equations \eqref{qn} and
\eqref{rn},
we get the equations
\begin{equation}
f=1+g+z(tf^2+t^2f^2) \quad \text{and} \quad g=zf+tzfg+t^2zfg.
\end{equation}
We then have that $G=zf$ satisfies the algebraic equation
\begin{equation}
(t+t^2)^2G^3-2(t+t^2)G^2+(1-z+z(t+t^2))G-z=0.
\end{equation}
\section{Closed Forms}
\label{closedformssection}
The best one can hope to achieve with respect to an enumerative problem is to
find a closed-form representation for that what is being enumerated. The main
tool we will use is the Lagrange-B\"{u}rmann Inversion Theorem. We will only use
a simplified version and state it here for completeness.
\begin{theo}[Lagrange-B\"{u}rmann Inversion Theorem]
\label{LBinversion}
Let $\phi(u)$ be a formal power series with $\phi_0 \neq 0$, and let $Y(z)$ be
the unique formal power series solution of the equation $Y=z\phi(Y)$. Then the
coefficient of $Y(z)$ of order $n$ is given by
\begin{equation}
[z^n]Y(z)=\frac{1}{n}[u^{n-1}]\phi(u)^n.
\end{equation}
\end{theo}
The notation $[z^n]\cdots$ means the coefficent of $z^n$ in the power series
expansion of the expression that
follows.
\subsection{The Simplified Problem}
We first consider this problem for the single-index sequence $c_d$.\par
To apply the Lagrange-B\"{u}rmann Inversion Theorem, we need to write the algebraic equation \eqref{algeq1} in the form
$G=z\phi(G)$ for some formal power series $\phi$. By simple arithmetic, one
arrives at $G=z\phi(G)$ with
\begin{equation}
\phi(G)=\frac{1-G}{(1-2G)^2}
\end{equation}
(note that $\phi(0)\neq 0$). Then by Theorem \ref{LBinversion},
\begin{equation}
\label{cdzvsu}
c_d=[z^d]G(z)=\frac{1}{d}[u^{d-1}]\left( \frac{1-u}{(1-2u)^2}  \right)^d.
\end{equation}
Since
\begin{equation}
(1-u)^d=\sum \limits_{n=0}^d\binom{d}{n}(-u)^n
\end{equation}
and
\begin{equation}
\frac{1}{(1-2u)^{2d}}=\sum
\limits_{n=0}^{\infty}\binom{n+2d-1}{n}(2u)^n,
\end{equation}
then Cauchy multiplication gives
\begin{align}
[u^{d-1}]\left( \frac{1-u}{(1-2u)^2}  \right)^d=&\sum \limits_{n=0}^{d-1}\binom{n+2d-1}{n}2^n \binom{d}{d-1-n}(-1)^{d-1-n} \nonumber \\
=&(-1)^{d-1}\sum \limits_{n=0}^{d-1}\frac{(n+2d-1)!}{n!(2d-1)!}\frac{d!(-1)^n2^n}{(d-1-n)!(n+1)!}\nonumber \\
=&(-1)^{d-1}\sum \limits_{n=0}^{d-1}\frac{(n+2d-1)!}{(2d-1)!}\frac{d!(-1)^n}{(d-1-n)!}\frac{2^n}{(2)_nn!}
\end{align}
since $(2)_n=2 \cdot 3 \cdots (2+n-1)=(n+1)!$, and $(-1)^n(1-d)_n=\frac{(d-1)!}{(d-n-1)!}$ gives
\begin{align}
=&d(-1)^{d-1}\sum \limits_{n=0}^{d-1}\frac{(n+2d-1)!}{(2d-1)!}\frac{(1-d)_n2^n}{(2)_nn!}\nonumber \\
=&d(-1)^{d-1}\sum \limits_{n=0}^{d-1}\frac{(2d)_n(1-d)_n}{(2)_n}\frac{2^n}{n!}\nonumber \\
=&d(-1)^{d-1}{}_2F_1([1-d,2d];[2];2).
\end{align}
We utilize \eqref{cdzvsu} and Euler's formula:
\begin{equation}
\label{Eulersform}
{}_2F_1([a,b];[c];z)=(1-z)^{-a} {}_2F_1\left([a,c-b];[c];\frac{z}{z-1}\right)
\end{equation}
(found in Gardner, et. al. \cite{GKP}, for instance) for the Gaussian hypergeometric function
\begin{equation}
{}_2F_1([a,b];[c];z)=\sum_{n=0}^{\infty}\frac{\left(a\right)_n\left(b\right)_n }{\left( c\right)_n}\frac{z^n}{n!}
\end{equation}
where $(x)_n=x(x+1)\cdots (x+n-1)$ is the Pochhammer symbol. With $z=2$, $c=2$, $a=1-d$, and $b=2-2d$, \eqref{cdzvsu} and \eqref{Eulersform} give the hypergeometric function stated in the following:
\begin{theo}
The number of combinatorial invariants for $\Xi_d$ is
\begin{equation}
c_d={}_2F_1([2-2d,1-d];[2];2).
\end{equation}
\end{theo}
This series has a
finite number of terms due to $1-d\in \mathbb{Z}^-$ being non-positive for $d \geq 1$.
\subsection{The Complete Problem}
\label{closedformcomplete}
We first find $[z^d]G(z,t)$. Following the same method as before, we note
\begin{equation}
\phi(G)=\frac{1+(1-(t+t^2))G}{(1-(t+t^2)G)^2},
\end{equation}
and by Theorem \ref{LBinversion}, we again have
\begin{equation}
c_d(t):=[z^d]G(z,t)=\frac{1}{d}[u^{d-1}]\left(\frac{1+(1-(t+t^2))u}{(1-(t+t^2)u)^2}
\right)^d.
\end{equation}
Again utilizing Cauchy multiplication, we have
\begin{equation}
c_d(t)=\frac{1}{d}\sum \limits_{n=0}^{d-1}\binom{n+2d-1}{n}(t+t^2)^n
\binom{d}{d-1-n}(1-t-t^2)^{d-1-n},
\end{equation}
and rewriting gives
\begin{equation}
=\sum\limits_{n=0}^{d-1}\frac{(-1)^n(2d)_n(1-d)_n}{(2)_nn!}(t+t^2)^n(1-t-t^2)^{d-1-n}.
\end{equation}
Recognizing that this is the right-hand side of \eqref{Eulersform}, with $a=1-d$, $b=2-2d$, $c=2$, and $z=t+t^2$, this simplifies beautifully to
\begin{equation}
c_d(t)={}_2F_1([1-d,2(1-d)];[2];t+t^2).
\end{equation}
By expanding $(t+t^2)^n$, we get
\begin{theo}
The number of combinatorial invariants of dimension $q$ for $\Xi_d$ is
\begin{equation}
c_{d,q}=[t^q]{}_2F_1([1-d,2(1-d)];[2];t+t^2)=\sum_{n=0}^{d-1}\frac{(2-2d)_n(1-d)_n}{(2)_n n!}\binom{n}{q-n}.
\end{equation}
\end{theo}
\section{Asymptotic Behavior}
\label{asymptoticbehaviorsection}
The next characteristic one might want to examine is the asymptotic behavior of
the sequences $c_d$ and $c_{d,q}$.
\subsection{The Simplified Problem}
\label{asymsimple}
For $c_d$, we wish to utilize the relation
\begin{equation}
\label{Rrelation}
\frac{1}{R}=\limsup_{d \rightarrow \infty}\sqrt[d]{|c_d|},
\end{equation}
where $R$ is the radius of convergence of the associated generating function. We
first check that $R>0$ by use of the Implicit Function Theorem. We have
\begin{align}
P(G,z)&=G^3-G^2+\left(\frac{z+1}{4}  \right)G-\frac{z}{4}=0 \\
\frac{\partial P}{\partial G}\restrictto{(0,0)}&=\frac{1}{4}\neq 0,
\end{align}
so it follows that $R>0$. We next determine $R$ exactly in order to make a more
precise asymptotic estimation for $c_d$. The polynomial $P(G)$ behaves nicely;
it has three roots, counted with multiplicity. That is, $G$ has three branches,
of which we are interested in the one going through $(z,G)=(0,0)$.  The radius of
convergence of the power series expansion of $G$ about $z=0$ will be determined
by the first singularity on the branch of interest. In this light, we
wish to find all of the places the branches intersect, i.e. where the
discriminant is zero. Calculation gives
\begin{equation}
\Delta(P(G))=\frac{-z}{16}\left(z- \left( \frac{11\pm 5\sqrt{5}}{-2}\right)  \right),
\end{equation}
so the discriminant $\Delta(P(G))=0$ when $z=0$ and $z=\frac{11\pm 5\sqrt{5}}{-2}$.
When $z=0$, $P(G)$ has a double root at $G=1/2$, and a single root at $G=0$. The
branch we are interested in is the one where $G=0$. So $G$ is analytic in a neighborhood of $z=0$, and the radius of convergence will be determined by the $z$ with smallest modulus where there is a singularity. There are two numbers to check; we start with the
one with smallest modulus, $\frac{11-5\sqrt{5}}{-2}$. Using Sturm's Theorem, it can be shown that $P(G)$
has three distinct real roots in the interval $]0,\frac{11-5\sqrt{5}}{-2}[$. At
$z=\frac{11-5\sqrt{5}}{-2}$, $P(G)$ has a single root at
$G=\frac{1}{2}(-1+\sqrt{5})$ and a double root at $G=\frac{1}{4}(3-\sqrt{5})$.
Since $\frac{1}{4}(3-\sqrt{5})<\frac{1}{2}(-1+\sqrt{5})$, we can conclude that
the branch we are interested in intersects the point $(z,G)=(\frac{11-5\sqrt{5}}{-2},\frac{1}{4}(3-\sqrt{5}))$. Hence,
$R=\frac{11-5\sqrt{5}}{-2}$. Now since $c_d$ is positive and increasing,
equation \eqref{Rrelation} becomes
\begin{equation}
\frac{1}{R}=\lim_{d \rightarrow \infty}\sqrt[d]{c_d},
\end{equation}
and by the definition of a limit, we can conclude
\begin{equation}
c_d \sim (1/R)^d =\left(\frac{-2}{11-5\sqrt{5}}  \right)^d\approx (11.09)^d.
\end{equation}
\subsection{The Complete Problem}
We consider the discrete probability distribution
\begin{equation}
\pi_{d,q}=\frac{c_{d,q}}{c_d}.
\end{equation}
The goal is to show that this probability
distribution is asymptotically of some known form, for example normal, when $d
\rightarrow \infty$. We wish to use the same method as in \cite{FN1999}; compare with their Theorem 5.\par
\begin{theo}
The distribution $\pi_{d,q}$ has mean $\mu_d$ and variance $\sigma^2_d$ that satisfy
\begin{equation}
\mu_d \sim \kappa d, \quad \sigma^2_d \sim \lambda d,
\end{equation}
where $\kappa$ and $\lambda$ are algebraic numbers. The laws in each case are
asymptotically normal for $d\rightarrow \infty$.
\end{theo}
\proc{Proof.}
The method of the proof is taken from \cite{FN1999} and also draws on the references
\cite{FS_analcomb}, \cite{EB1973}, and \cite{FO1990}.\par
Consider the characteristic function of $\pi_d(q)$, normalized by shifting the mean
by $\mu_d$ and the variance by $\sigma_d^2$:
\begin{equation}
 f_d(t)= e^{-i\mu_dt/\sigma_d}G_d\left(e^{it/\sigma_d} \right),
\end{equation}
where $G_d(t)=\frac{c_d(t)}{c_d}$ is the probability generating function.
The goal is to show that the normalized characteristic functions $f_d(t)$ converge
pointwise to the characteristic function $e^{-t^2/2}$ of the standard normal. Then
the limit law follows by Levy's Continuity Theorem.  \par
We look at
\begin{equation}
 P(z,t,G)=(t+t^2)^2G^3-2(t+t^2)G^2+(1-z+z(t+t^2))G-z.
\end{equation}
By previous discussions, $G_0=(3-\sqrt{5})/4$ is a double root for $P(\rho,1,G)$, where $\rho = \frac{11-5 \sqrt{5}}{-2}$, so
by Weierstrass preparation, there is an analytic factorization
\begin{equation}
P(z,t,G)=((G-G_0)^2+m_1(z,t)(G-G_0)+m_2(z,t))\cdot H(z,t,G),
\end{equation}
where $H(z,t,G)$ is analytic near $(\rho,1,G_0)$, $H(\rho,1,G_0)\neq 0$, and $m_1$ and $m_2$
are analytic at $(\rho,1)$.
Then applying the quadratic formula,
\begin{equation}
G(z,t)-G_0=\frac{1}{2}\left( -m_1(z,t)\pm\sqrt{m_1(z,t)^2-4m_2(z,t)}  \right).
\end{equation}
We determine our branch of interest (as we did in Subsection \ref{strstablecase}) by the following. Consider first $(z,t)$
restricted to $0\leq z<\rho(t)$ and $0\leq t < 1$. Since $G(z,t)$ is real there ($G$
has real coefficients), the discriminant $m_1(z,t)^2-4m_2(z,t)$ must also be real
and non-negative. Furthermore, since $G$ is increasing in $z$ for fixed $t$
(coefficients of $G$ are non-negative), and since the discriminant vanishes at
$\rho(t)$ and is hence decreasing in $z$, then we need to take the minus sign, i.e.
 \begin{equation}
G(z,t)-G_0=\frac{1}{2}\left( -m_1(z,t)-\sqrt{m_1(z,t)^2-4m_2(z,t)}  \right).
\end{equation}
Set $D(z,t):=m_1(z,t)^2-4m_2(z,t)$. Consider the resultant
\begin{equation}
R(z,t)={\rm Result}_G \left( P(z,t,G), \frac{\partial}{\partial G}P(z,t,G) \right),
\end{equation}
which is a polynomial whose restriction $R(z,1)$ has, by the discussion in Section
\ref{asymsimple}, a simple root at $\rho = \frac{11-5 \sqrt{5}}{-2}$.  By the
Implicit Function Theorem, this root lifts to an analytic branch $\rho(t)$ of an
algebraic function, for $t$ near 1:
\begin{equation}
R(\rho(t),t)=0, \quad \rho(1)=\rho.
\end{equation}
Therefore, the function $D(z,1)$ has a simple real zero at $z=\rho$, and we can
factorize
\begin{equation}
 D(z,t)=(\rho(t)-z)K(z,t),
\end{equation}
for some analytic $K$ satisfying $K(\rho,1)\neq0$. So the uniform family of singular
expansions
for $G(z,t)$ near $(\rho(t),t)$ takes the form
\begin{equation}
\label{puiexp}
G(z,t)=c_0(t)+c_1(t)\sqrt{1-z/\rho(t)}+\mathcal{O}(1-z/\rho(t)),
\end{equation}
uniformly with respect to $t$ for $t$ in a small neighborhood of 1, and with
$\rho(t)$, $c_0(t)$, and $c_1(t)$ analytic at 1. If $G(z,t)$ satisfies equation
\eqref{puiexp}, we let $\tilde{z}=z/\rho(t)$ and apply the results in \cite{FO1990}
to obtain
\begin{align}
 c_d(t):=[z^d]G(z,t)&=(\rho(t))^{-d}[\tilde{z}^d]\{
c_1(t)\sqrt{1-\tilde{z}}+\mathcal{O}(1-\tilde{z})\}\nonumber \\
&=\frac{c_1(t)(\rho(t))^{-d}}{\Gamma(-1/2)d^{3/2}}\left( 1 + \mathcal{O}(1/d)\right).
\end{align}
 \par
The probability generating function is therefore a so-called ``quasi-power,'' i.e.
$G_d(t)$ satisfies
\begin{equation}
\label{quasipower}
G_d(t)=\frac{\gamma(t)}{\gamma(1)}\left( \frac{\rho(1)}{\rho(t)} \right)^d\left(
1+\mathcal{O}\left( \frac{1}{d} \right) \right).
\end{equation}
 If \eqref{quasipower} holds, then for $d \rightarrow \infty$ and close enough to
$w=1$,
\begin{equation}
 G_d(t)\sim \left(\frac{\rho(1)}{\rho(t)}\right)^d,
\end{equation}
so for fixed $t$, letting $d\rightarrow \infty$,
\begin{equation}
 f_d(t) \sim e^{-i\mu_dt/\sigma_d}\left[ \frac{r(0)}{r(it/\sigma_d)}\right]^d,
\end{equation}
where $r(s)=\rho(e^s)$.
 Set $y(s)=\log (r(s)/r(0))$, expand $y(s)$ in a Taylor series, and use $\left[
\frac{r(0)}{r(s)}\right]^d=\exp \left(d \log (r(0)/r(s))\right)$ to get
\begin{equation}
 f_d(t)\sim \exp\{-i\mu_dt/\sigma_d -d \left(ity'(0)/\sigma_d
-t^2y''(0)/2\sigma_d^2+\mathcal{O}(t^3/\sigma_d^3) \right) \}.
\end{equation}
By calculation on \eqref{quasipower},
\begin{align}
 \mu_d &= G_d'(1)\sim -d \rho'(1)/\rho(1)=-dy'(0)\quad \text{and}\nonumber \\
\sigma_d^2&=G_d''(1)+\mu_d(1-\mu_d)\sim -d
\left(\frac{\rho''(1)+\rho'(1)}{\rho(1)}-\left(
\frac{\rho'(1)}{\rho(1)}\right)^2\right)=-dy''(0).
\end{align}
We obtain
\begin{equation}
 f_d(t)\sim e^{-t^2/2}
\end{equation}
for all $t$.
\ep \medbreak
By differentiation of \eqref{quasipower} and with help from MAPLE's
\emph{algeqtodiffeq} function in the GFUN package \cite{SZ1994}, we can calculate
\begin{equation}
\rho(1) = \frac{-11+5\sqrt{5}}{2}, \quad \rho'(1)=\frac{87-39\sqrt{5}}{2}, \quad
\rho''(1)=\frac{702\rho(1)+716\rho'(1)}{-60},
\end{equation}
so
\begin{equation}
\kappa=\frac{\rho'(1)}{\rho(1)}=\frac{3(-3+\sqrt{5})}{-2}\approx 1.145898036,
\end{equation}
and
\begin{equation}
\lambda=-\frac{\rho''(1)}{\rho(1)}-\frac{\rho'(1)}{\rho(1)}+\left(
\frac{\rho'(1)}{\rho(1)}\right)^2=\frac{-60+29\sqrt{5}}{10}\approx 0.484597133.
\end{equation}
\section{Alternative Counting Problems}
\label{alternativeprobssection}
There are several variations on this enumeration problem that one might consider
interesting. We discuss a few of these now.
\subsection{Distinguishing between Real and Complex Dimension}
It has been suggested that one might want to count the complex dimensions
(corresponding to the number of analytic invariants in $\mathbb{H}$ for a vector field $\xi_P$) and the real
dimensions (corresponding to the number of analytic invariants in $\mathbb{R}_+$ for a vector field $\xi_P$)
separately (see \cite{BD09} and \cite{Sent} for definitions of the analytic invariants).  This would correspond in the complete problem to a triple-indexed
sequence $c_{d,s,h}$, where $s$ is the number of complex analytic invariants and $h$
is the number of real analytic invariants. We follow the same method as in Section
\ref{genfctsection}.
Let
\begin{equation}
s_{n,s,h}=\begin{cases}
"s"s_{n-1,s,h}&\text{if string does not start with parentheses } \\
"(s"s_{2a,s_1,h_1} "s)" s_{b,s_2,h_2}& s_1+s_2=s, \ h_1+h_2+1=h\\
"[s"s_{2a,s_1,h_1} "s]" s_{b,s_2,h_2}& s_1+s_2+1=s, \ h_1+h_2=h\end{cases},
\end{equation}
be our string generating algorithm for $2a+b+2=n$. We define a recursion equation
for $p_{n,s,h}$ in the usual way,
and arrive at the generating function $G(z,t,w)=\sum c_{d,s,h}z^d w^s t^h$
satisfying (unsurprisingly)
\begin{equation}
 (t+w)^2G^3-2(t+w)G^2+(1-z+z(t+w))G-z=0.
\end{equation}
One can use methods as before to determine properties of $c_{d,s,h}$.
\subsection{Enumeration of Moduli Classes}
The combinatorial classes we have focused on distinguish between vector fields which
are conformally conjugate by a rotation by a $(d-1)^{{\rm st}}$ root of unity. Some
may consider the problem of enumeration of classes in moduli space (i.e. without the
labelling of the separatrices) more valuable. In order to study this, we utilize
P\'{o}lya theory for the cyclic group isomorphic to the group generated by $e^{2\pi
{\rm i}/(d-1)}$. We are interested in the sequence $ \tilde{c}_d=[z^{d-1}]C(G(z))$
where $C(G(z))$ is generating function for the number of combinatorial classes for
degree $d$ vector fields moduli rotation by a $(d-1)^{{\rm st}}$ root of unity.
P\'{o}lya theory gives
\begin{equation}
 C(G(z))=\sum \limits_{k \geq 1}\frac{\phi(k)}{k}\log \frac{1}{1-G(z^k)},
\end{equation}
where $\phi(k)$ is the Euler totient function. This relies on prime factorization,
so it does not seem likely that one can find a closed-form expression for
$\tilde{c}_d$ as a function of $d$ as of date.

\end{document}

%% file: trajsatinfty.pstex_t
\begin{picture}(0,0)%
\includegraphics{trajsatinfty.pstex}%
\end{picture}%
\setlength{\unitlength}{3947sp}%
\begingroup\makeatletter\ifx\SetFigFont\undefined%
\gdef\SetFigFont#1#2#3#4#5{%
  \reset@font\fontsize{#1}{#2pt}%
  \fontfamily{#3}\fontseries{#4}\fontshape{#5}%
  \selectfont}%
\fi\endgroup%
\begin{picture}(4755,4902)(511,-4057)
\put(5251,-1636){\makebox(0,0)[lb]{\smash{{\SetFigFont{17}{20.4}{\rmdefault}{\bfdefault}{\updefault}{\color[rgb]{0,0,0}$\gamma_0$}%
}}}}
\put(3151,-1711){\makebox(0,0)[lb]{\smash{{\SetFigFont{17}{20.4}{\rmdefault}{\bfdefault}{\updefault}{\color[rgb]{0,0,0}$\infty$}%
}}}}
\put(526,-1636){\makebox(0,0)[lb]{\smash{{\SetFigFont{17}{20.4}{\rmdefault}{\bfdefault}{\updefault}{\color[rgb]{0,0,0}$\gamma_4$}%
}}}}
\put(4576,-3211){\makebox(0,0)[lb]{\smash{{\SetFigFont{17}{20.4}{\rmdefault}{\bfdefault}{\updefault}{\color[rgb]{0,0,0}$\gamma_1$}%
}}}}
\put(4576,-61){\makebox(0,0)[lb]{\smash{{\SetFigFont{17}{20.4}{\rmdefault}{\bfdefault}{\updefault}{\color[rgb]{0,0,0}$\gamma_7$}%
}}}}
\put(2926,-3961){\makebox(0,0)[lb]{\smash{{\SetFigFont{17}{20.4}{\rmdefault}{\bfdefault}{\updefault}{\color[rgb]{0,0,0}$\gamma_2$}%
}}}}
\put(3001,614){\makebox(0,0)[lb]{\smash{{\SetFigFont{17}{20.4}{\rmdefault}{\bfdefault}{\updefault}{\color[rgb]{0,0,0}$\gamma_6$}%
}}}}
\put(1276,-3361){\makebox(0,0)[lb]{\smash{{\SetFigFont{17}{20.4}{\rmdefault}{\bfdefault}{\updefault}{\color[rgb]{0,0,0}$\gamma_3$}%
}}}}
\put(1276, 89){\makebox(0,0)[lb]{\smash{{\SetFigFont{17}{20.4}{\rmdefault}{\bfdefault}{\updefault}{\color[rgb]{0,0,0}$\gamma_5$}%
}}}}
\end{picture}%

%% file: trajsatinfty2.pstex_t
\begin{picture}(0,0)%
\includegraphics{trajsatinfty2.pstex}%
\end{picture}%
\setlength{\unitlength}{3947sp}%
\begingroup\makeatletter\ifx\SetFigFont\undefined%
\gdef\SetFigFont#1#2#3#4#5{%
  \reset@font\fontsize{#1}{#2pt}%
  \fontfamily{#3}\fontseries{#4}\fontshape{#5}%
  \selectfont}%
\fi\endgroup%
\begin{picture}(4755,4902)(511,-4057)
\put(5251,-1636){\makebox(0,0)[lb]{\smash{{\SetFigFont{17}{20.4}{\rmdefault}{\bfdefault}{\updefault}{\color[rgb]{0,0,0}$\gamma_0$}%
}}}}
\put(3151,-1711){\makebox(0,0)[lb]{\smash{{\SetFigFont{17}{20.4}{\rmdefault}{\bfdefault}{\updefault}{\color[rgb]{0,0,0}$\infty$}%
}}}}
\put(526,-1636){\makebox(0,0)[lb]{\smash{{\SetFigFont{17}{20.4}{\rmdefault}{\bfdefault}{\updefault}{\color[rgb]{0,0,0}$\gamma_4$}%
}}}}
\put(4576, 14){\makebox(0,0)[lb]{\smash{{\SetFigFont{17}{20.4}{\rmdefault}{\bfdefault}{\updefault}{\color[rgb]{0,0,0}$\gamma_7$}%
}}}}
\put(4651,-3286){\makebox(0,0)[lb]{\smash{{\SetFigFont{17}{20.4}{\rmdefault}{\bfdefault}{\updefault}{\color[rgb]{0,0,0}$\gamma_1$}%
}}}}
\put(2926,-3961){\makebox(0,0)[lb]{\smash{{\SetFigFont{17}{20.4}{\rmdefault}{\bfdefault}{\updefault}{\color[rgb]{0,0,0}$\gamma_2$}%
}}}}
\put(2851,614){\makebox(0,0)[lb]{\smash{{\SetFigFont{17}{20.4}{\rmdefault}{\bfdefault}{\updefault}{\color[rgb]{0,0,0}$\gamma_6$}%
}}}}
\put(1351, 89){\makebox(0,0)[lb]{\smash{{\SetFigFont{17}{20.4}{\rmdefault}{\bfdefault}{\updefault}{\color[rgb]{0,0,0}$\gamma_5$}%
}}}}
\put(1276,-3361){\makebox(0,0)[lb]{\smash{{\SetFigFont{17}{20.4}{\rmdefault}{\bfdefault}{\updefault}{\color[rgb]{0,0,0}$\gamma_3$}%
}}}}
\put(3526,-1936){\makebox(0,0)[lb]{\smash{{\SetFigFont{17}{20.4}{\rmdefault}{\bfdefault}{\updefault}{\color[rgb]{0,0,0}$e_1$}%
}}}}
\put(3451,-1411){\makebox(0,0)[lb]{\smash{{\SetFigFont{17}{20.4}{\rmdefault}{\bfdefault}{\updefault}{\color[rgb]{0,0,0}$e_0$}%
}}}}
\put(2626,-1036){\makebox(0,0)[lb]{\smash{{\SetFigFont{17}{20.4}{\rmdefault}{\bfdefault}{\updefault}{\color[rgb]{0,0,0}$e_6$}%
}}}}
\put(2251,-1411){\makebox(0,0)[lb]{\smash{{\SetFigFont{17}{20.4}{\rmdefault}{\bfdefault}{\updefault}{\color[rgb]{0,0,0}$e_5$}%
}}}}
\put(2251,-1861){\makebox(0,0)[lb]{\smash{{\SetFigFont{17}{20.4}{\rmdefault}{\bfdefault}{\updefault}{\color[rgb]{0,0,0}$e_4$}%
}}}}
\put(2626,-2161){\makebox(0,0)[lb]{\smash{{\SetFigFont{17}{20.4}{\rmdefault}{\bfdefault}{\updefault}{\color[rgb]{0,0,0}$e_3$}%
}}}}
\put(3076,-2161){\makebox(0,0)[lb]{\smash{{\SetFigFont{17}{20.4}{\rmdefault}{\bfdefault}{\updefault}{\color[rgb]{0,0,0}$e_2$}%
}}}}
\put(3076,-1111){\makebox(0,0)[lb]{\smash{{\SetFigFont{17}{20.4}{\rmdefault}{\bfdefault}{\updefault}{\color[rgb]{0,0,0}$e_7$}%
}}}}
\end{picture}%

%% file: separatrixgraph2.pstex_t
\begin{picture}(0,0)%
\includegraphics{separatrixgraph2.pstex}%
\end{picture}%
\setlength{\unitlength}{3947sp}%
\begingroup\makeatletter\ifx\SetFigFont\undefined%
\gdef\SetFigFont#1#2#3#4#5{%
  \reset@font\fontsize{#1}{#2pt}%
  \fontfamily{#3}\fontseries{#4}\fontshape{#5}%
  \selectfont}%
\fi\endgroup%
\begin{picture}(7080,7468)(586,-6563)
\put(1501,-5386){\makebox(0,0)[lb]{\smash{{\SetFigFont{20}{24.0}{\rmdefault}{\mddefault}{\updefault}{\color[rgb]{0,0,0}$s_5$}%
}}}}
\put(601,-2836){\makebox(0,0)[lb]{\smash{{\SetFigFont{20}{24.0}{\rmdefault}{\mddefault}{\updefault}{\color[rgb]{0,0,0}$s_4$}%
}}}}
\put(4126,614){\makebox(0,0)[lb]{\smash{{\SetFigFont{20}{24.0}{\rmdefault}{\mddefault}{\updefault}{\color[rgb]{0,0,0}$s_2$}%
}}}}
\put(6676,-436){\makebox(0,0)[lb]{\smash{{\SetFigFont{20}{24.0}{\rmdefault}{\mddefault}{\updefault}{\color[rgb]{0,0,0}$s_1$}%
}}}}
\put(7651,-2836){\makebox(0,0)[lb]{\smash{{\SetFigFont{20}{24.0}{\rmdefault}{\mddefault}{\updefault}{\color[rgb]{0,0,0}$s_{3,0}$}%
}}}}
\put(4126,-6436){\makebox(0,0)[lb]{\smash{{\SetFigFont{20}{24.0}{\rmdefault}{\mddefault}{\updefault}{\color[rgb]{0,0,0}$s_{7,6}$}%
}}}}
\end{picture}%

%% file: centerzoneshade1.pstex_t
\begin{picture}(0,0)%
\includegraphics{centerzoneshade1.pstex}%
\end{picture}%
\setlength{\unitlength}{3947sp}%
\begingroup\makeatletter\ifx\SetFigFont\undefined%
\gdef\SetFigFont#1#2#3#4#5{%
  \reset@font\fontsize{#1}{#2pt}%
  \fontfamily{#3}\fontseries{#4}\fontshape{#5}%
  \selectfont}%
\fi\endgroup%
\begin{picture}(7500,7500)(526,-7261)
\put(7351,-3211){\makebox(0,0)[lb]{\smash{{\SetFigFont{20}{24.0}{\rmdefault}{\mddefault}{\updefault}{\color[rgb]{0,0,0}$s_{5,0}$}%
}}}}
\put(7351,-1486){\makebox(0,0)[lb]{\smash{{\SetFigFont{20}{24.0}{\rmdefault}{\mddefault}{\updefault}{\color[rgb]{0,0,0}$e_1$}%
}}}}
\put(4426,-3361){\makebox(0,0)[lb]{\smash{{\SetFigFont{20}{24.0}{\rmdefault}{\mddefault}{\updefault}{\color[rgb]{0,0,0}$\zeta=$center}%
}}}}
\put(1276,-61){\makebox(0,0)[lb]{\smash{{\SetFigFont{20}{24.0}{\rmdefault}{\mddefault}{\updefault}{\color[rgb]{0,0,0}$s_{1,2}$}%
}}}}
\put(601,-2161){\makebox(0,0)[lb]{\smash{{\SetFigFont{20}{24.0}{\rmdefault}{\mddefault}{\updefault}{\color[rgb]{0,0,0}$e_3$}%
}}}}
\put(2626,-7111){\makebox(0,0)[lb]{\smash{{\SetFigFont{20}{24.0}{\rmdefault}{\mddefault}{\updefault}{\color[rgb]{0,0,0}$s_{3,4}$}%
}}}}
\put(4501,-7111){\makebox(0,0)[lb]{\smash{{\SetFigFont{20}{24.0}{\rmdefault}{\mddefault}{\updefault}{\color[rgb]{0,0,0}$e_5$}%
}}}}
\end{picture}%

%% file: sepalzoneshade2.pstex_t
\begin{picture}(0,0)%
\includegraphics{sepalzoneshade2.pstex}%
\end{picture}%
\setlength{\unitlength}{3947sp}%
\begingroup\makeatletter\ifx\SetFigFont\undefined%
\gdef\SetFigFont#1#2#3#4#5{%
  \reset@font\fontsize{#1}{#2pt}%
  \fontfamily{#3}\fontseries{#4}\fontshape{#5}%
  \selectfont}%
\fi\endgroup%
\begin{picture}(7931,7588)(214,-7305)
\put(4276,-61){\makebox(0,0)[lb]{\smash{{\SetFigFont{20}{24.0}{\rmdefault}{\mddefault}{\updefault}{\color[rgb]{0,0,0}$s_{3,2}$}%
}}}}
\put(5701,-61){\makebox(0,0)[lb]{\smash{{\SetFigFont{20}{24.0}{\rmdefault}{\mddefault}{\updefault}{\color[rgb]{0,0,0}$e_2$}%
}}}}
\put(7276,-3886){\makebox(0,0)[lb]{\smash{{\SetFigFont{20}{24.0}{\rmdefault}{\mddefault}{\updefault}{\color[rgb]{0,0,0}$s_{1,0}$}%
}}}}
\put(7501,-5161){\makebox(0,0)[lb]{\smash{{\SetFigFont{20}{24.0}{\rmdefault}{\mddefault}{\updefault}{\color[rgb]{0,0,0}$e_0$}%
}}}}
\put(1426,-7111){\makebox(0,0)[lb]{\smash{{\SetFigFont{20}{24.0}{\rmdefault}{\mddefault}{\updefault}{\color[rgb]{0,0,0}$s_5$}%
}}}}
\put(3751,-7111){\makebox(0,0)[lb]{\smash{{\SetFigFont{20}{24.0}{\rmdefault}{\mddefault}{\updefault}{\color[rgb]{0,0,0}$s_{7,6}$}%
}}}}
\put(751,-3286){\makebox(0,0)[lb]{\smash{{\SetFigFont{20}{24.0}{\rmdefault}{\mddefault}{\updefault}{\color[rgb]{0,0,0}$s_4$}%
}}}}
\put(2476,-7111){\makebox(0,0)[lb]{\smash{{\SetFigFont{20}{24.0}{\rmdefault}{\mddefault}{\updefault}{\color[rgb]{0,0,0}$e_6$}%
}}}}
\put(4051,-2836){\makebox(0,0)[lb]{\smash{{\SetFigFont{20}{24.0}{\rmdefault}{\mddefault}{\updefault}{\color[rgb]{0,0,0}limit point}%
}}}}
\put(3301,-2836){\makebox(0,0)[lb]{\smash{{\SetFigFont{20}{24.0}{\rmdefault}{\mddefault}{\updefault}{\color[rgb]{0,0,0}$\alpha, \omega$}%
}}}}
\put(751,-2011){\makebox(0,0)[lb]{\smash{{\SetFigFont{20}{24.0}{\rmdefault}{\mddefault}{\updefault}{\color[rgb]{0,0,0}$e_4$}%
}}}}
\end{picture}%

%% file: alphaomegashade2.pstex_t
\begin{picture}(0,0)%
\includegraphics{alphaomegashade2.pstex}%
\end{picture}%
\setlength{\unitlength}{3947sp}%
\begingroup\makeatletter\ifx\SetFigFont\undefined%
\gdef\SetFigFont#1#2#3#4#5{%
  \reset@font\fontsize{#1}{#2pt}%
  \fontfamily{#3}\fontseries{#4}\fontshape{#5}%
  \selectfont}%
\fi\endgroup%
\begin{picture}(7544,7544)(407,-7261)
\put(3451,-7111){\makebox(0,0)[lb]{\smash{{\SetFigFont{20}{24.0}{\rmdefault}{\mddefault}{\updefault}{\color[rgb]{0,0,0}$s_{7,6}$}%
}}}}
\put(7351,-3211){\makebox(0,0)[lb]{\smash{{\SetFigFont{20}{24.0}{\rmdefault}{\mddefault}{\updefault}{\color[rgb]{0,0,0}$s_0$}%
}}}}
\put(7201,-811){\makebox(0,0)[lb]{\smash{{\SetFigFont{20}{24.0}{\rmdefault}{\mddefault}{\updefault}{\color[rgb]{0,0,0}$s_{1,2}$}%
}}}}
\put(526,-736){\makebox(0,0)[lb]{\smash{{\SetFigFont{20}{24.0}{\rmdefault}{\mddefault}{\updefault}{\color[rgb]{0,0,0}$s_3$}%
}}}}
\put(6826,-4936){\makebox(0,0)[lb]{\smash{{\SetFigFont{20}{24.0}{\rmdefault}{\mddefault}{\updefault}{\color[rgb]{0,0,0}$e_0$}%
}}}}
\put(2401,-136){\makebox(0,0)[lb]{\smash{{\SetFigFont{20}{24.0}{\rmdefault}{\mddefault}{\updefault}{\color[rgb]{0,0,0}$e_3$}%
}}}}
\put(6976,-1786){\makebox(0,0)[lb]{\smash{{\SetFigFont{20}{24.0}{\rmdefault}{\mddefault}{\updefault}{\color[rgb]{0,0,0}$e_1$}%
}}}}
\put(5026,-2161){\makebox(0,0)[lb]{\smash{{\SetFigFont{20}{24.0}{\rmdefault}{\mddefault}{\updefault}{\color[rgb]{0,0,0}$\alpha$ limit point}%
}}}}
\put(601,-2611){\makebox(0,0)[lb]{\smash{{\SetFigFont{20}{24.0}{\rmdefault}{\mddefault}{\updefault}{\color[rgb]{0,0,0}$e_4$}%
}}}}
\put(901,-7111){\makebox(0,0)[lb]{\smash{{\SetFigFont{20}{24.0}{\rmdefault}{\mddefault}{\updefault}{\color[rgb]{0,0,0}$s_{5,4}$}%
}}}}
\put(2101,-7036){\makebox(0,0)[lb]{\smash{{\SetFigFont{20}{24.0}{\rmdefault}{\mddefault}{\updefault}{\color[rgb]{0,0,0}$e_6$}%
}}}}
\put(2776,-4561){\makebox(0,0)[lb]{\smash{{\SetFigFont{20}{24.0}{\rmdefault}{\mddefault}{\updefault}{\color[rgb]{0,0,0}$\omega$ limit point}%
}}}}
\end{picture}%

%% file: alphaomegashade.pstex_t
\begin{picture}(0,0)%
\includegraphics{alphaomegashade.pstex}%
\end{picture}%
\setlength{\unitlength}{3947sp}%
\begingroup\makeatletter\ifx\SetFigFont\undefined%
\gdef\SetFigFont#1#2#3#4#5{%
  \reset@font\fontsize{#1}{#2pt}%
  \fontfamily{#3}\fontseries{#4}\fontshape{#5}%
  \selectfont}%
\fi\endgroup%
\begin{picture}(14627,7544)(407,-7261)
\put(12301,-3211){\makebox(0,0)[lb]{\smash{{\SetFigFont{20}{24.0}{\rmdefault}{\mddefault}{\updefault}{\color[rgb]{0,0,0}$\frac{d}{dz}$}%
}}}}
\put(10801,-2386){\makebox(0,0)[lb]{\smash{{\SetFigFont{20}{24.0}{\rmdefault}{\mddefault}{\updefault}{\color[rgb]{0,0,0}$S_0$}%
}}}}
\put(10801,-4561){\makebox(0,0)[lb]{\smash{{\SetFigFont{20}{24.0}{\rmdefault}{\mddefault}{\updefault}{\color[rgb]{0,0,0}$S_0$}%
}}}}
\put(14551,-2386){\makebox(0,0)[lb]{\smash{{\SetFigFont{20}{24.0}{\rmdefault}{\mddefault}{\updefault}{\color[rgb]{0,0,0}$S_3$}%
}}}}
\put(14551,-4561){\makebox(0,0)[lb]{\smash{{\SetFigFont{20}{24.0}{\rmdefault}{\mddefault}{\updefault}{\color[rgb]{0,0,0}$S_3$}%
}}}}
\put(11851,-2836){\makebox(0,0)[lb]{\smash{{\SetFigFont{17}{20.4}{\rmdefault}{\mddefault}{\updefault}{\color[rgb]{0,0,0}$E_0$}%
}}}}
\put(11926,-4111){\makebox(0,0)[lb]{\smash{{\SetFigFont{17}{20.4}{\rmdefault}{\mddefault}{\updefault}{\color[rgb]{0,0,0}$E_1$}%
}}}}
\put(12826,-4111){\makebox(0,0)[lb]{\smash{{\SetFigFont{17}{20.4}{\rmdefault}{\mddefault}{\updefault}{\color[rgb]{0,0,0}$E_3$}%
}}}}
\put(12826,-2836){\makebox(0,0)[lb]{\smash{{\SetFigFont{17}{20.4}{\rmdefault}{\mddefault}{\updefault}{\color[rgb]{0,0,0}$E_6$}%
}}}}
\put(13726,-2836){\makebox(0,0)[lb]{\smash{{\SetFigFont{17}{20.4}{\rmdefault}{\mddefault}{\updefault}{\color[rgb]{0,0,0}$E_4$}%
}}}}
\put(12226,-2386){\makebox(0,0)[lb]{\smash{{\SetFigFont{17}{20.4}{\rmdefault}{\mddefault}{\updefault}{\color[rgb]{0,0,0}$S_{7,6}$}%
}}}}
\put(12151,-4561){\makebox(0,0)[lb]{\smash{{\SetFigFont{17}{20.4}{\rmdefault}{\mddefault}{\updefault}{\color[rgb]{0,0,0}$S_{1,2}$}%
}}}}
\put(13126,-2386){\makebox(0,0)[lb]{\smash{{\SetFigFont{17}{20.4}{\rmdefault}{\mddefault}{\updefault}{\color[rgb]{0,0,0}$S_{5,4}$}%
}}}}
\put(6976,-2161){\makebox(0,0)[lb]{\smash{{\SetFigFont{20}{24.0}{\rmdefault}{\mddefault}{\updefault}{\color[rgb]{0,0,0}$e_1$}%
}}}}
\put(7351,-3211){\makebox(0,0)[lb]{\smash{{\SetFigFont{20}{24.0}{\rmdefault}{\mddefault}{\updefault}{\color[rgb]{0,0,0}$s_0$}%
}}}}
\put(976,-7111){\makebox(0,0)[lb]{\smash{{\SetFigFont{20}{24.0}{\rmdefault}{\mddefault}{\updefault}{\color[rgb]{0,0,0}$s_{5,4}$}%
}}}}
\put(2101,-7036){\makebox(0,0)[lb]{\smash{{\SetFigFont{20}{24.0}{\rmdefault}{\mddefault}{\updefault}{\color[rgb]{0,0,0}$e_6$}%
}}}}
\put(6976,-6211){\makebox(0,0)[lb]{\smash{{\SetFigFont{20}{24.0}{\rmdefault}{\mddefault}{\updefault}{\color[rgb]{0,0,0}$s_{7,6}$}%
}}}}
\put(6826,-4936){\makebox(0,0)[lb]{\smash{{\SetFigFont{20}{24.0}{\rmdefault}{\mddefault}{\updefault}{\color[rgb]{0,0,0}$e_0$}%
}}}}
\put(9001,-3136){\makebox(0,0)[lb]{\smash{{\SetFigFont{20}{24.0}{\rmdefault}{\mddefault}{\updefault}{\color[rgb]{0,0,0}$\phi_{e_1}$}%
}}}}
\put(3901,-3436){\makebox(0,0)[lb]{\smash{{\SetFigFont{20}{24.0}{\rmdefault}{\mddefault}{\updefault}{\color[rgb]{0,0,0}$P(z)\frac{d}{dz}$}%
}}}}
\put(601,-2611){\makebox(0,0)[lb]{\smash{{\SetFigFont{20}{24.0}{\rmdefault}{\mddefault}{\updefault}{\color[rgb]{0,0,0}$e_4$}%
}}}}
\put(526,-736){\makebox(0,0)[lb]{\smash{{\SetFigFont{20}{24.0}{\rmdefault}{\mddefault}{\updefault}{\color[rgb]{0,0,0}$s_3$}%
}}}}
\put(2401,-136){\makebox(0,0)[lb]{\smash{{\SetFigFont{20}{24.0}{\rmdefault}{\mddefault}{\updefault}{\color[rgb]{0,0,0}$e_3$}%
}}}}
\put(7276,-886){\makebox(0,0)[lb]{\smash{{\SetFigFont{20}{24.0}{\rmdefault}{\mddefault}{\updefault}{\color[rgb]{0,0,0}$s_{1,2}$}%
}}}}
\end{picture}%

%% file: sepalzoneshade.pstex_t
\begin{picture}(0,0)%
\includegraphics{sepalzoneshade.pstex}%
\end{picture}%
\setlength{\unitlength}{3947sp}%
\begingroup\makeatletter\ifx\SetFigFont\undefined%
\gdef\SetFigFont#1#2#3#4#5{%
  \reset@font\fontsize{#1}{#2pt}%
  \fontfamily{#3}\fontseries{#4}\fontshape{#5}%
  \selectfont}%
\fi\endgroup%
\begin{picture}(14177,7588)(557,-7305)
\put(10426,-2161){\makebox(0,0)[lb]{\smash{{\SetFigFont{20}{24.0}{\rmdefault}{\mddefault}{\updefault}{\color[rgb]{0,0,0}$S_4$}%
}}}}
\put(14701,-2161){\makebox(0,0)[lb]{\smash{{\SetFigFont{20}{24.0}{\rmdefault}{\mddefault}{\updefault}{\color[rgb]{0,0,0}$S_5$}%
}}}}
\put(11176,-2536){\makebox(0,0)[lb]{\smash{{\SetFigFont{20}{24.0}{\rmdefault}{\mddefault}{\updefault}{\color[rgb]{0,0,0}$E_4$}%
}}}}
\put(12076,-2536){\makebox(0,0)[lb]{\smash{{\SetFigFont{20}{24.0}{\rmdefault}{\mddefault}{\updefault}{\color[rgb]{0,0,0}$E_2$}%
}}}}
\put(13126,-2536){\makebox(0,0)[lb]{\smash{{\SetFigFont{20}{24.0}{\rmdefault}{\mddefault}{\updefault}{\color[rgb]{0,0,0}$E_0$}%
}}}}
\put(14176,-2536){\makebox(0,0)[lb]{\smash{{\SetFigFont{20}{24.0}{\rmdefault}{\mddefault}{\updefault}{\color[rgb]{0,0,0}$E_6$}%
}}}}
\put(11401,-2011){\makebox(0,0)[lb]{\smash{{\SetFigFont{20}{24.0}{\rmdefault}{\mddefault}{\updefault}{\color[rgb]{0,0,0}$S_{3,2}$}%
}}}}
\put(12451,-2011){\makebox(0,0)[lb]{\smash{{\SetFigFont{20}{24.0}{\rmdefault}{\mddefault}{\updefault}{\color[rgb]{0,0,0}$S_{1,0}$}%
}}}}
\put(13576,-2011){\makebox(0,0)[lb]{\smash{{\SetFigFont{20}{24.0}{\rmdefault}{\mddefault}{\updefault}{\color[rgb]{0,0,0}$S_{7,6}$}%
}}}}
\put(1426,-7111){\makebox(0,0)[lb]{\smash{{\SetFigFont{20}{24.0}{\rmdefault}{\mddefault}{\updefault}{\color[rgb]{0,0,0}$s_5$}%
}}}}
\put(2476,-7111){\makebox(0,0)[lb]{\smash{{\SetFigFont{20}{24.0}{\rmdefault}{\mddefault}{\updefault}{\color[rgb]{0,0,0}$e_6$}%
}}}}
\put(3751,-7111){\makebox(0,0)[lb]{\smash{{\SetFigFont{20}{24.0}{\rmdefault}{\mddefault}{\updefault}{\color[rgb]{0,0,0}$s_{7,6}$}%
}}}}
\put(7501,-5161){\makebox(0,0)[lb]{\smash{{\SetFigFont{20}{24.0}{\rmdefault}{\mddefault}{\updefault}{\color[rgb]{0,0,0}$e_0$}%
}}}}
\put(7276,-3886){\makebox(0,0)[lb]{\smash{{\SetFigFont{20}{24.0}{\rmdefault}{\mddefault}{\updefault}{\color[rgb]{0,0,0}$s_{1,0}$}%
}}}}
\put(9001,-3136){\makebox(0,0)[lb]{\smash{{\SetFigFont{20}{24.0}{\rmdefault}{\mddefault}{\updefault}{\color[rgb]{0,0,0}$\phi_{e_4}$}%
}}}}
\put(3601,-2761){\makebox(0,0)[lb]{\smash{{\SetFigFont{20}{24.0}{\rmdefault}{\mddefault}{\updefault}{\color[rgb]{0,0,0}$P(z)\frac{d}{dz}$}%
}}}}
\put(751,-2011){\makebox(0,0)[lb]{\smash{{\SetFigFont{20}{24.0}{\rmdefault}{\mddefault}{\updefault}{\color[rgb]{0,0,0}$e_4$}%
}}}}
\put(751,-3286){\makebox(0,0)[lb]{\smash{{\SetFigFont{20}{24.0}{\rmdefault}{\mddefault}{\updefault}{\color[rgb]{0,0,0}$s_4$}%
}}}}
\put(5701,-61){\makebox(0,0)[lb]{\smash{{\SetFigFont{20}{24.0}{\rmdefault}{\mddefault}{\updefault}{\color[rgb]{0,0,0}$e_2$}%
}}}}
\put(4276,-61){\makebox(0,0)[lb]{\smash{{\SetFigFont{20}{24.0}{\rmdefault}{\mddefault}{\updefault}{\color[rgb]{0,0,0}$s_{3,2}$}%
}}}}
\end{picture}%

%% file: centerzoneshade.pstex_t
\begin{picture}(0,0)%
\includegraphics{centerzoneshade.pstex}%
\end{picture}%
\setlength{\unitlength}{3947sp}%
\begingroup\makeatletter\ifx\SetFigFont\undefined%
\gdef\SetFigFont#1#2#3#4#5{%
  \reset@font\fontsize{#1}{#2pt}%
  \fontfamily{#3}\fontseries{#4}\fontshape{#5}%
  \selectfont}%
\fi\endgroup%
\begin{picture}(13905,7500)(586,-7261)
\put(10876,-4111){\makebox(0,0)[lb]{\smash{{\SetFigFont{17}{20.4}{\rmdefault}{\mddefault}{\updefault}{\color[rgb]{0,0,0}$E_1$}%
}}}}
\put(12076,-4111){\makebox(0,0)[lb]{\smash{{\SetFigFont{17}{20.4}{\rmdefault}{\mddefault}{\updefault}{\color[rgb]{0,0,0}$E_3$}%
}}}}
\put(13276,-4111){\makebox(0,0)[lb]{\smash{{\SetFigFont{17}{20.4}{\rmdefault}{\mddefault}{\updefault}{\color[rgb]{0,0,0}$E_5$}%
}}}}
\put(14476,-4111){\makebox(0,0)[lb]{\smash{{\SetFigFont{17}{20.4}{\rmdefault}{\mddefault}{\updefault}{\color[rgb]{0,0,0}$E_1$}%
}}}}
\put(11026,-4561){\makebox(0,0)[lb]{\smash{{\SetFigFont{17}{20.4}{\rmdefault}{\mddefault}{\updefault}{\color[rgb]{0,0,0}$S_{1,2}$}%
}}}}
\put(12301,-4561){\makebox(0,0)[lb]{\smash{{\SetFigFont{20}{24.0}{\rmdefault}{\mddefault}{\updefault}{\color[rgb]{0,0,0}$S_{3,4}$}%
}}}}
\put(13651,-4561){\makebox(0,0)[lb]{\smash{{\SetFigFont{20}{24.0}{\rmdefault}{\mddefault}{\updefault}{\color[rgb]{0,0,0}$S_{5,0}$}%
}}}}
\put(12076,-2836){\makebox(0,0)[lb]{\smash{{\SetFigFont{20}{24.0}{\rmdefault}{\mddefault}{\updefault}{\color[rgb]{0,0,0}$\frac{d}{dz}$}%
}}}}
\put(2626,-7111){\makebox(0,0)[lb]{\smash{{\SetFigFont{20}{24.0}{\rmdefault}{\mddefault}{\updefault}{\color[rgb]{0,0,0}$s_{3,4}$}%
}}}}
\put(4501,-7111){\makebox(0,0)[lb]{\smash{{\SetFigFont{20}{24.0}{\rmdefault}{\mddefault}{\updefault}{\color[rgb]{0,0,0}$e_5$}%
}}}}
\put(9001,-3136){\makebox(0,0)[lb]{\smash{{\SetFigFont{20}{24.0}{\rmdefault}{\mddefault}{\updefault}{\color[rgb]{0,0,0}$\phi_{e_1}$}%
}}}}
\put(7351,-3211){\makebox(0,0)[lb]{\smash{{\SetFigFont{20}{24.0}{\rmdefault}{\mddefault}{\updefault}{\color[rgb]{0,0,0}$s_{5,0}$}%
}}}}
\put(3601,-2761){\makebox(0,0)[lb]{\smash{{\SetFigFont{20}{24.0}{\rmdefault}{\mddefault}{\updefault}{\color[rgb]{0,0,0}$P(z)\frac{d}{dz}$}%
}}}}
\put(7351,-1486){\makebox(0,0)[lb]{\smash{{\SetFigFont{20}{24.0}{\rmdefault}{\mddefault}{\updefault}{\color[rgb]{0,0,0}$e_1$}%
}}}}
\put(1276,-61){\makebox(0,0)[lb]{\smash{{\SetFigFont{20}{24.0}{\rmdefault}{\mddefault}{\updefault}{\color[rgb]{0,0,0}$s_{1,2}$}%
}}}}
\put(601,-2161){\makebox(0,0)[lb]{\smash{{\SetFigFont{20}{24.0}{\rmdefault}{\mddefault}{\updefault}{\color[rgb]{0,0,0}$e_3$}%
}}}}
\end{picture}%

%% file: transversals2_defensepres.pstex_t
\begin{picture}(0,0)%
\includegraphics{transversals2_defensepres.pstex}%
\end{picture}%
\setlength{\unitlength}{3947sp}%
\begingroup\makeatletter\ifx\SetFigFont\undefined%
\gdef\SetFigFont#1#2#3#4#5{%
  \reset@font\fontsize{#1}{#2pt}%
  \fontfamily{#3}\fontseries{#4}\fontshape{#5}%
  \selectfont}%
\fi\endgroup%
\begin{picture}(7577,7566)(407,-7283)
\put(526,-736){\makebox(0,0)[lb]{\smash{{\SetFigFont{20}{24.0}{\rmdefault}{\mddefault}{\updefault}{\color[rgb]{0,0,0}$s_3$}%
}}}}
\put(7351,-3211){\makebox(0,0)[lb]{\smash{{\SetFigFont{20}{24.0}{\rmdefault}{\mddefault}{\updefault}{\color[rgb]{0,0,0}$s_0$}%
}}}}
\put(7276,-586){\makebox(0,0)[lb]{\smash{{\SetFigFont{20}{24.0}{\rmdefault}{\mddefault}{\updefault}{\color[rgb]{0,0,0}$s_{1,2}$}%
}}}}
\put(3451,-7111){\makebox(0,0)[lb]{\smash{{\SetFigFont{20}{24.0}{\rmdefault}{\mddefault}{\updefault}{\color[rgb]{0,0,0}$s_{7,6}$}%
}}}}
\put(526,-3586){\makebox(0,0)[lb]{\smash{{\SetFigFont{20}{24.0}{\rmdefault}{\mddefault}{\updefault}{\color[rgb]{0,0,0}$s_{5,4}$}%
}}}}
\end{picture}%

%% file: rectmap.pstex_t
\begin{picture}(0,0)%
\includegraphics{rectmap.pstex}%
\end{picture}%
\setlength{\unitlength}{3947sp}%
\begingroup\makeatletter\ifx\SetFigFont\undefined%
\gdef\SetFigFont#1#2#3#4#5{%
  \reset@font\fontsize{#1}{#2pt}%
  \fontfamily{#3}\fontseries{#4}\fontshape{#5}%
  \selectfont}%
\fi\endgroup%
\begin{picture}(937,740)(1861,-1620)
\put(1876,-1111){\makebox(0,0)[lb]{\smash{{\SetFigFont{17}{20.4}{\rmdefault}{\mddefault}{\updefault}{\color[rgb]{0,0,0}$\int_{e_1}^{z}\frac{\rm d \zeta}{P(\zeta)} $}%
}}}}
\end{picture}%

%% file: transversals_diststrip.pstex_t
\begin{picture}(0,0)%
\includegraphics{transversals_diststrip.pstex}%
\end{picture}%
\setlength{\unitlength}{3947sp}%
\begingroup\makeatletter\ifx\SetFigFont\undefined%
\gdef\SetFigFont#1#2#3#4#5{%
  \reset@font\fontsize{#1}{#2pt}%
  \fontfamily{#3}\fontseries{#4}\fontshape{#5}%
  \selectfont}%
\fi\endgroup%
\begin{picture}(3930,1803)(1261,-2155)
\put(3451,-1336){\makebox(0,0)[lb]{\smash{{\SetFigFont{14}{16.8}{\rmdefault}{\mddefault}{\updefault}{\color[rgb]{1,0,1}$T_{3,0}$}%
}}}}
\put(3976,-1786){\makebox(0,0)[lb]{\smash{{\SetFigFont{12}{14.4}{\rmdefault}{\mddefault}{\updefault}{\color[rgb]{0,0,0}$e_3$}%
}}}}
\put(5176,-1936){\makebox(0,0)[lb]{\smash{{\SetFigFont{12}{14.4}{\rmdefault}{\mddefault}{\updefault}{\color[rgb]{0,0,0}$s_3$}%
}}}}
\put(3526,-886){\makebox(0,0)[lb]{\smash{{\SetFigFont{12}{14.4}{\rmdefault}{\mddefault}{\updefault}{\color[rgb]{0,0,0}$e_6$}%
}}}}
\put(4426,-886){\makebox(0,0)[lb]{\smash{{\SetFigFont{12}{14.4}{\rmdefault}{\mddefault}{\updefault}{\color[rgb]{0,0,0}$e_4$}%
}}}}
\put(5176,-736){\makebox(0,0)[lb]{\smash{{\SetFigFont{12}{14.4}{\rmdefault}{\mddefault}{\updefault}{\color[rgb]{0,0,0}$s_3$}%
}}}}
\put(2026,-1411){\makebox(0,0)[lb]{\smash{{\SetFigFont{12}{14.4}{\rmdefault}{\mddefault}{\updefault}{\color[rgb]{0,0,0}$\frac{d}{dz}$}%
}}}}
\put(2626,-1786){\makebox(0,0)[lb]{\smash{{\SetFigFont{12}{14.4}{\rmdefault}{\mddefault}{\updefault}{\color[rgb]{0,0,0}$e_1$}%
}}}}
\put(2326,-886){\makebox(0,0)[lb]{\smash{{\SetFigFont{12}{14.4}{\rmdefault}{\mddefault}{\updefault}{\color[rgb]{0,0,0}$e_0$}%
}}}}
\put(1276,-736){\makebox(0,0)[lb]{\smash{{\SetFigFont{12}{14.4}{\rmdefault}{\mddefault}{\updefault}{\color[rgb]{0,0,0}$s_0$}%
}}}}
\put(1276,-1936){\makebox(0,0)[lb]{\smash{{\SetFigFont{12}{14.4}{\rmdefault}{\mddefault}{\updefault}{\color[rgb]{0,0,0}$s_0$}%
}}}}
\put(3151,-2086){\makebox(0,0)[lb]{\smash{{\SetFigFont{12}{14.4}{\rmdefault}{\mddefault}{\updefault}{\color[rgb]{0,0,0}$s_{1,2}$}%
}}}}
\put(2851,-511){\makebox(0,0)[lb]{\smash{{\SetFigFont{12}{14.4}{\rmdefault}{\mddefault}{\updefault}{\color[rgb]{0,0,0}$s_{7,6}$}%
}}}}
\put(3901,-511){\makebox(0,0)[lb]{\smash{{\SetFigFont{12}{14.4}{\rmdefault}{\mddefault}{\updefault}{\color[rgb]{0,0,0}$s_{5,4}$}%
}}}}
\end{picture}%

%% file: disttransvf.pstex_t
\begin{picture}(0,0)%
\includegraphics{disttransvf.pstex}%
\end{picture}%
\setlength{\unitlength}{3947sp}%
\begingroup\makeatletter\ifx\SetFigFont\undefined%
\gdef\SetFigFont#1#2#3#4#5{%
  \reset@font\fontsize{#1}{#2pt}%
  \fontfamily{#3}\fontseries{#4}\fontshape{#5}%
  \selectfont}%
\fi\endgroup%
\begin{picture}(10473,8359)(61,-7754)
\put(526,-736){\makebox(0,0)[lb]{\smash{{\SetFigFont{20}{24.0}{\rmdefault}{\mddefault}{\updefault}{\color[rgb]{0,0,0}$s_3$}%
}}}}
\put(7351,-3211){\makebox(0,0)[lb]{\smash{{\SetFigFont{20}{24.0}{\rmdefault}{\mddefault}{\updefault}{\color[rgb]{0,0,0}$s_0$}%
}}}}
\put(601,-3511){\makebox(0,0)[lb]{\smash{{\SetFigFont{20}{24.0}{\rmdefault}{\mddefault}{\updefault}{\color[rgb]{0,0,0}$s_{5,4}$}%
}}}}
\put(2926,-7036){\makebox(0,0)[lb]{\smash{{\SetFigFont{20}{24.0}{\rmdefault}{\mddefault}{\updefault}{\color[rgb]{0,0,0}$s_{7,6}$}%
}}}}
\put(6976,-961){\makebox(0,0)[lb]{\smash{{\SetFigFont{20}{24.0}{\rmdefault}{\mddefault}{\updefault}{\color[rgb]{0,0,0}$s_{1,2}$}%
}}}}
\put(8026,-5236){\makebox(0,0)[lb]{\smash{{\SetFigFont{20}{24.0}{\rmdefault}{\mddefault}{\updefault}{\color[rgb]{0,0,0}$e_0$}%
}}}}
\put(8026,-2086){\makebox(0,0)[lb]{\smash{{\SetFigFont{20}{24.0}{\rmdefault}{\mddefault}{\updefault}{\color[rgb]{0,0,0}$e_1$}%
}}}}
\put(2551,314){\makebox(0,0)[lb]{\smash{{\SetFigFont{20}{24.0}{\rmdefault}{\mddefault}{\updefault}{\color[rgb]{0,0,0}$e_3$}%
}}}}
\put( 76,-2311){\makebox(0,0)[lb]{\smash{{\SetFigFont{20}{24.0}{\rmdefault}{\mddefault}{\updefault}{\color[rgb]{0,0,0}$e_4$}%
}}}}
\put(2026,-7636){\makebox(0,0)[lb]{\smash{{\SetFigFont{20}{24.0}{\rmdefault}{\mddefault}{\updefault}{\color[rgb]{0,0,0}$e_6$}%
}}}}
\put(4726,-2686){\makebox(0,0)[lb]{\smash{{\SetFigFont{20}{24.0}{\rmdefault}{\mddefault}{\updefault}{\color[rgb]{.69,0,.69}$T_{3,0}$}%
}}}}
\put(3451,-4111){\makebox(0,0)[lb]{\smash{{\SetFigFont{20}{24.0}{\rmdefault}{\mddefault}{\updefault}{\color[rgb]{0,0,0}$P(z)\frac{d}{dz}$}%
}}}}
\put(9001,-3136){\makebox(0,0)[lb]{\smash{{\SetFigFont{25}{30.0}{\rmdefault}{\mddefault}{\updefault}{\color[rgb]{0,0,0}$\int_{\infty}^z \frac{d\zeta}{P(\zeta)}$}%
}}}}
\end{picture}%

%% file: disttrans.pstex_t
\begin{picture}(0,0)%
\includegraphics{disttrans.pstex}%
\end{picture}%
\setlength{\unitlength}{3947sp}%
\begingroup\makeatletter\ifx\SetFigFont\undefined%
\gdef\SetFigFont#1#2#3#4#5{%
  \reset@font\fontsize{#1}{#2pt}%
  \fontfamily{#3}\fontseries{#4}\fontshape{#5}%
  \selectfont}%
\fi\endgroup%
\begin{picture}(3930,2860)(1261,-3212)
\put(3451,-1336){\makebox(0,0)[lb]{\smash{{\SetFigFont{14}{16.8}{\rmdefault}{\mddefault}{\updefault}{\color[rgb]{1,0,1}$T_{3,0}$}%
}}}}
\put(3976,-1786){\makebox(0,0)[lb]{\smash{{\SetFigFont{12}{14.4}{\rmdefault}{\mddefault}{\updefault}{\color[rgb]{0,0,0}$e_3$}%
}}}}
\put(5176,-1936){\makebox(0,0)[lb]{\smash{{\SetFigFont{12}{14.4}{\rmdefault}{\mddefault}{\updefault}{\color[rgb]{0,0,0}$s_3$}%
}}}}
\put(3526,-886){\makebox(0,0)[lb]{\smash{{\SetFigFont{12}{14.4}{\rmdefault}{\mddefault}{\updefault}{\color[rgb]{0,0,0}$e_6$}%
}}}}
\put(4426,-886){\makebox(0,0)[lb]{\smash{{\SetFigFont{12}{14.4}{\rmdefault}{\mddefault}{\updefault}{\color[rgb]{0,0,0}$e_4$}%
}}}}
\put(5176,-736){\makebox(0,0)[lb]{\smash{{\SetFigFont{12}{14.4}{\rmdefault}{\mddefault}{\updefault}{\color[rgb]{0,0,0}$s_3$}%
}}}}
\put(2026,-1411){\makebox(0,0)[lb]{\smash{{\SetFigFont{12}{14.4}{\rmdefault}{\mddefault}{\updefault}{\color[rgb]{0,0,0}$\frac{d}{dz}$}%
}}}}
\put(2626,-1786){\makebox(0,0)[lb]{\smash{{\SetFigFont{12}{14.4}{\rmdefault}{\mddefault}{\updefault}{\color[rgb]{0,0,0}$e_1$}%
}}}}
\put(2326,-886){\makebox(0,0)[lb]{\smash{{\SetFigFont{12}{14.4}{\rmdefault}{\mddefault}{\updefault}{\color[rgb]{0,0,0}$e_0$}%
}}}}
\put(1276,-736){\makebox(0,0)[lb]{\smash{{\SetFigFont{12}{14.4}{\rmdefault}{\mddefault}{\updefault}{\color[rgb]{0,0,0}$s_0$}%
}}}}
\put(1276,-1936){\makebox(0,0)[lb]{\smash{{\SetFigFont{12}{14.4}{\rmdefault}{\mddefault}{\updefault}{\color[rgb]{0,0,0}$s_0$}%
}}}}
\put(3151,-2086){\makebox(0,0)[lb]{\smash{{\SetFigFont{12}{14.4}{\rmdefault}{\mddefault}{\updefault}{\color[rgb]{0,0,0}$s_{1,2}$}%
}}}}
\put(2851,-511){\makebox(0,0)[lb]{\smash{{\SetFigFont{12}{14.4}{\rmdefault}{\mddefault}{\updefault}{\color[rgb]{0,0,0}$s_{7,6}$}%
}}}}
\put(3901,-511){\makebox(0,0)[lb]{\smash{{\SetFigFont{12}{14.4}{\rmdefault}{\mddefault}{\updefault}{\color[rgb]{0,0,0}$s_{5,4}$}%
}}}}
\end{picture}%

%% file: transversalgraph.pstex_t
\begin{picture}(0,0)%
\includegraphics{transversalgraph.pstex}%
\end{picture}%
\setlength{\unitlength}{3947sp}%
\begingroup\makeatletter\ifx\SetFigFont\undefined%
\gdef\SetFigFont#1#2#3#4#5{%
  \reset@font\fontsize{#1}{#2pt}%
  \fontfamily{#3}\fontseries{#4}\fontshape{#5}%
  \selectfont}%
\fi\endgroup%
\begin{picture}(7080,7525)(586,-6872)
\put(3151,-2836){\makebox(0,0)[lb]{\smash{{\SetFigFont{20}{24.0}{\rmdefault}{\mddefault}{\updefault}{\color[rgb]{0,0,0}$T_{5,4}$}%
}}}}
\put(676,-1636){\makebox(0,0)[lb]{\smash{{\SetFigFont{20}{24.0}{\rmdefault}{\mddefault}{\updefault}{\color[rgb]{0,0,0}$e_4$}%
}}}}
\put(601,-4561){\makebox(0,0)[lb]{\smash{{\SetFigFont{20}{24.0}{\rmdefault}{\mddefault}{\updefault}{\color[rgb]{0,0,0}$e_5$}%
}}}}
\put(2476,-6511){\makebox(0,0)[lb]{\smash{{\SetFigFont{20}{24.0}{\rmdefault}{\mddefault}{\updefault}{\color[rgb]{0,0,0}$e_6$}%
}}}}
\put(5401,-6436){\makebox(0,0)[lb]{\smash{{\SetFigFont{20}{24.0}{\rmdefault}{\mddefault}{\updefault}{\color[rgb]{0,0,0}$e_7$}%
}}}}
\put(7201,-4711){\makebox(0,0)[lb]{\smash{{\SetFigFont{20}{24.0}{\rmdefault}{\mddefault}{\updefault}{\color[rgb]{0,0,0}$e_0$}%
}}}}
\put(1501,-5686){\makebox(0,0)[lb]{\smash{{\SetFigFont{20}{24.0}{\rmdefault}{\mddefault}{\updefault}{\color[rgb]{0,0,0}$s_5$}%
}}}}
\put(601,-3136){\makebox(0,0)[lb]{\smash{{\SetFigFont{20}{24.0}{\rmdefault}{\mddefault}{\updefault}{\color[rgb]{0,0,0}$s_4$}%
}}}}
\put(4126,314){\makebox(0,0)[lb]{\smash{{\SetFigFont{20}{24.0}{\rmdefault}{\mddefault}{\updefault}{\color[rgb]{0,0,0}$s_2$}%
}}}}
\put(6676,-736){\makebox(0,0)[lb]{\smash{{\SetFigFont{20}{24.0}{\rmdefault}{\mddefault}{\updefault}{\color[rgb]{0,0,0}$s_1$}%
}}}}
\put(7651,-3136){\makebox(0,0)[lb]{\smash{{\SetFigFont{20}{24.0}{\rmdefault}{\mddefault}{\updefault}{\color[rgb]{0,0,0}$s_{3,0}$}%
}}}}
\put(4126,-6736){\makebox(0,0)[lb]{\smash{{\SetFigFont{20}{24.0}{\rmdefault}{\mddefault}{\updefault}{\color[rgb]{0,0,0}$s_{7,6}$}%
}}}}
\put(7426,-1861){\makebox(0,0)[lb]{\smash{{\SetFigFont{20}{24.0}{\rmdefault}{\mddefault}{\updefault}{\color[rgb]{0,0,0}$e_1$}%
}}}}
\put(5776,-61){\makebox(0,0)[lb]{\smash{{\SetFigFont{20}{24.0}{\rmdefault}{\mddefault}{\updefault}{\color[rgb]{0,0,0}$e_2$}%
}}}}
\put(2626,164){\makebox(0,0)[lb]{\smash{{\SetFigFont{20}{24.0}{\rmdefault}{\mddefault}{\updefault}{\color[rgb]{0,0,0}$e_3$}%
}}}}
\end{picture}%

%% file: openHchains1.pstex_t
\begin{picture}(0,0)%
\includegraphics{openHchains1.pstex}%
\end{picture}%
\setlength{\unitlength}{3947sp}%
\begingroup\makeatletter\ifx\SetFigFont\undefined%
\gdef\SetFigFont#1#2#3#4#5{%
  \reset@font\fontsize{#1}{#2pt}%
  \fontfamily{#3}\fontseries{#4}\fontshape{#5}%
  \selectfont}%
\fi\endgroup%
\begin{picture}(12423,3459)(661,-2644)
\put(5851,314){\makebox(0,0)[lb]{\smash{{\SetFigFont{20}{24.0}{\rmdefault}{\mddefault}{\updefault}{\color[rgb]{0,0,0}$s_{k_2,j_2}$}%
}}}}
\put(9451,-61){\makebox(0,0)[lb]{\smash{{\SetFigFont{20}{24.0}{\rmdefault}{\mddefault}{\updefault}{\color[rgb]{0,0,0}$s_{k_1,j_1}$}%
}}}}
\put(12751,464){\makebox(0,0)[lb]{\smash{{\SetFigFont{20}{24.0}{\rmdefault}{\mddefault}{\updefault}{\color[rgb]{0,0,0}$s_{j_1-1}$}%
}}}}
\put(676, 89){\makebox(0,0)[lb]{\smash{{\SetFigFont{20}{24.0}{\rmdefault}{\mddefault}{\updefault}{\color[rgb]{0,0,0}$s_{k_3+1}$}%
}}}}
\put(2701,164){\makebox(0,0)[lb]{\smash{{\SetFigFont{20}{24.0}{\rmdefault}{\mddefault}{\updefault}{\color[rgb]{0,0,0}$s_{k_3,j_3}$}%
}}}}
\end{picture}%

%% file: openHchains2.pstex_t
\begin{picture}(0,0)%
\includegraphics{openHchains2.pstex}%
\end{picture}%
\setlength{\unitlength}{3947sp}%
\begingroup\makeatletter\ifx\SetFigFont\undefined%
\gdef\SetFigFont#1#2#3#4#5{%
  \reset@font\fontsize{#1}{#2pt}%
  \fontfamily{#3}\fontseries{#4}\fontshape{#5}%
  \selectfont}%
\fi\endgroup%
\begin{picture}(12423,3447)(661,-2644)
\put(12751,464){\makebox(0,0)[lb]{\smash{{\SetFigFont{20}{24.0}{\rmdefault}{\mddefault}{\updefault}{\color[rgb]{0,0,0}$s_{k_1-1}$}%
}}}}
\put(9451,-61){\makebox(0,0)[lb]{\smash{{\SetFigFont{20}{24.0}{\rmdefault}{\mddefault}{\updefault}{\color[rgb]{0,0,0}$s_{k_1,j_1}$}%
}}}}
\put(2701,164){\makebox(0,0)[lb]{\smash{{\SetFigFont{20}{24.0}{\rmdefault}{\mddefault}{\updefault}{\color[rgb]{0,0,0}$s_{k_3,j_3}$}%
}}}}
\put(5851,314){\makebox(0,0)[lb]{\smash{{\SetFigFont{20}{24.0}{\rmdefault}{\mddefault}{\updefault}{\color[rgb]{0,0,0}$s_{k_2,j_2}$}%
}}}}
\put(676, 89){\makebox(0,0)[lb]{\smash{{\SetFigFont{20}{24.0}{\rmdefault}{\mddefault}{\updefault}{\color[rgb]{0,0,0}$s_{j_3+1}$}%
}}}}
\end{picture}%

%% file: openHchains.pstex_t
\begin{picture}(0,0)%
\includegraphics{openHchains.pstex}%
\end{picture}%
\setlength{\unitlength}{3947sp}%
\begingroup\makeatletter\ifx\SetFigFont\undefined%
\gdef\SetFigFont#1#2#3#4#5{%
  \reset@font\fontsize{#1}{#2pt}%
  \fontfamily{#3}\fontseries{#4}\fontshape{#5}%
  \selectfont}%
\fi\endgroup%
\begin{picture}(12423,7305)(661,-6502)
\put(10201,-6361){\makebox(0,0)[lb]{\smash{{\SetFigFont{20}{24.0}{\rmdefault}{\mddefault}{\updefault}{\color[rgb]{0,0,0}$s_j$}%
}}}}
\put(9451,-61){\makebox(0,0)[lb]{\smash{{\SetFigFont{20}{24.0}{\rmdefault}{\mddefault}{\updefault}{\color[rgb]{0,0,0}$s_{k_1,j_1}$}%
}}}}
\put(2701,164){\makebox(0,0)[lb]{\smash{{\SetFigFont{20}{24.0}{\rmdefault}{\mddefault}{\updefault}{\color[rgb]{0,0,0}$s_{k_3,j_3}$}%
}}}}
\put(5851,314){\makebox(0,0)[lb]{\smash{{\SetFigFont{20}{24.0}{\rmdefault}{\mddefault}{\updefault}{\color[rgb]{0,0,0}$s_{k_2,j_2}$}%
}}}}
\put(676, 89){\makebox(0,0)[lb]{\smash{{\SetFigFont{20}{24.0}{\rmdefault}{\mddefault}{\updefault}{\color[rgb]{0,0,0}$s_{j_3+1}$}%
}}}}
\put(12751,464){\makebox(0,0)[lb]{\smash{{\SetFigFont{20}{24.0}{\rmdefault}{\mddefault}{\updefault}{\color[rgb]{0,0,0}$s_{k_1-1}$}%
}}}}
\put(6226,-2161){\makebox(0,0)[lb]{\smash{{\SetFigFont{20}{24.0}{\rmdefault}{\mddefault}{\updefault}{\color[rgb]{0,0,0}$T_{j_3+1,j}$}%
}}}}
\end{picture}%

%% file: d4generic.pstex_t
\begin{picture}(0,0)%
\includegraphics{d4generic.pstex}%
\end{picture}%
\setlength{\unitlength}{3947sp}%
\begingroup\makeatletter\ifx\SetFigFont\undefined%
\gdef\SetFigFont#1#2#3#4#5{%
  \reset@font\fontsize{#1}{#2pt}%
  \fontfamily{#3}\fontseries{#4}\fontshape{#5}%
  \selectfont}%
\fi\endgroup%
\begin{picture}(10189,7431)(286,-6827)
\put(301,-996){\makebox(0,0)[lb]{\smash{{\SetFigFont{11}{13.2}{\rmdefault}{\bfdefault}{\updefault}{\color[rgb]{0,0,0}$s_3$}%
}}}}
\put(2999,-307){\makebox(0,0)[lb]{\smash{{\SetFigFont{11}{13.2}{\rmdefault}{\bfdefault}{\updefault}{\color[rgb]{0,0,0}$e_1$}%
}}}}
\put(2470,222){\makebox(0,0)[lb]{\smash{{\SetFigFont{11}{13.2}{\rmdefault}{\bfdefault}{\updefault}{\color[rgb]{0,0,0}$s_1$}%
}}}}
\put(1730,433){\makebox(0,0)[lb]{\smash{{\SetFigFont{11}{13.2}{\rmdefault}{\bfdefault}{\updefault}{\color[rgb]{0,0,0}$e_2$}%
}}}}
\put(990,222){\makebox(0,0)[lb]{\smash{{\SetFigFont{11}{13.2}{\rmdefault}{\bfdefault}{\updefault}{\color[rgb]{0,0,0}$s_2$}%
}}}}
\put(407,-307){\makebox(0,0)[lb]{\smash{{\SetFigFont{11}{13.2}{\rmdefault}{\bfdefault}{\updefault}{\color[rgb]{0,0,0}$e_3$}%
}}}}
\put(459,-1630){\makebox(0,0)[lb]{\smash{{\SetFigFont{11}{13.2}{\rmdefault}{\bfdefault}{\updefault}{\color[rgb]{0,0,0}$e_4$}%
}}}}
\put(1730,-2372){\makebox(0,0)[lb]{\smash{{\SetFigFont{11}{13.2}{\rmdefault}{\bfdefault}{\updefault}{\color[rgb]{0,0,0}$e_5$}%
}}}}
\put(2470,-2159){\makebox(0,0)[lb]{\smash{{\SetFigFont{11}{13.2}{\rmdefault}{\bfdefault}{\updefault}{\color[rgb]{0,0,0}$s_5$}%
}}}}
\put(3151,-1036){\makebox(0,0)[lb]{\smash{{\SetFigFont{11}{13.2}{\rmdefault}{\bfdefault}{\updefault}{\color[rgb]{0,0,0}$s_0$}%
}}}}
\put(10460,-1001){\makebox(0,0)[lb]{\smash{{\SetFigFont{11}{13.2}{\rmdefault}{\bfdefault}{\updefault}{\color[rgb]{0,0,0}$s_0$}%
}}}}
\put(10247,-313){\makebox(0,0)[lb]{\smash{{\SetFigFont{11}{13.2}{\rmdefault}{\bfdefault}{\updefault}{\color[rgb]{0,0,0}$e_1$}%
}}}}
\put(9720,218){\makebox(0,0)[lb]{\smash{{\SetFigFont{11}{13.2}{\rmdefault}{\bfdefault}{\updefault}{\color[rgb]{0,0,0}$s_1$}%
}}}}
\put(8979,428){\makebox(0,0)[lb]{\smash{{\SetFigFont{11}{13.2}{\rmdefault}{\bfdefault}{\updefault}{\color[rgb]{0,0,0}$e_2$}%
}}}}
\put(7656,-313){\makebox(0,0)[lb]{\smash{{\SetFigFont{11}{13.2}{\rmdefault}{\bfdefault}{\updefault}{\color[rgb]{0,0,0}$e_3$}%
}}}}
\put(7501,-1036){\makebox(0,0)[lb]{\smash{{\SetFigFont{11}{13.2}{\rmdefault}{\bfdefault}{\updefault}{\color[rgb]{0,0,0}$s_3$}%
}}}}
\put(1051,-2236){\makebox(0,0)[lb]{\smash{{\SetFigFont{11}{13.2}{\rmdefault}{\bfdefault}{\updefault}{\color[rgb]{0,0,0}$s_4$}%
}}}}
\put(3001,-1636){\makebox(0,0)[lb]{\smash{{\SetFigFont{11}{13.2}{\rmdefault}{\bfdefault}{\updefault}{\color[rgb]{0,0,0}$e_0$}%
}}}}
\put(6764,-1020){\makebox(0,0)[lb]{\smash{{\SetFigFont{11}{13.2}{\rmdefault}{\bfdefault}{\updefault}{\color[rgb]{0,0,0}$s_0$}%
}}}}
\put(5281,408){\makebox(0,0)[lb]{\smash{{\SetFigFont{11}{13.2}{\rmdefault}{\bfdefault}{\updefault}{\color[rgb]{0,0,0}$e_2$}%
}}}}
\put(4540,198){\makebox(0,0)[lb]{\smash{{\SetFigFont{11}{13.2}{\rmdefault}{\bfdefault}{\updefault}{\color[rgb]{0,0,0}$s_2$}%
}}}}
\put(4645,-2237){\makebox(0,0)[lb]{\smash{{\SetFigFont{11}{13.2}{\rmdefault}{\bfdefault}{\updefault}{\color[rgb]{0,0,0}$s_4$}%
}}}}
\put(5281,-2395){\makebox(0,0)[lb]{\smash{{\SetFigFont{11}{13.2}{\rmdefault}{\bfdefault}{\updefault}{\color[rgb]{0,0,0}$e_5$}%
}}}}
\put(3853,-1020){\makebox(0,0)[lb]{\smash{{\SetFigFont{11}{13.2}{\rmdefault}{\bfdefault}{\updefault}{\color[rgb]{0,0,0}$s_3$}%
}}}}
\put(4051,-286){\makebox(0,0)[lb]{\smash{{\SetFigFont{11}{13.2}{\rmdefault}{\bfdefault}{\updefault}{\color[rgb]{0,0,0}$e_3$}%
}}}}
\put(4051,-1636){\makebox(0,0)[lb]{\smash{{\SetFigFont{11}{13.2}{\rmdefault}{\bfdefault}{\updefault}{\color[rgb]{0,0,0}$e_4$}%
}}}}
\put(6076,-2236){\makebox(0,0)[lb]{\smash{{\SetFigFont{11}{13.2}{\rmdefault}{\bfdefault}{\updefault}{\color[rgb]{0,0,0}$s_5$}%
}}}}
\put(6601,-1711){\makebox(0,0)[lb]{\smash{{\SetFigFont{11}{13.2}{\rmdefault}{\bfdefault}{\updefault}{\color[rgb]{0,0,0}$e_0$}%
}}}}
\put(6676,-286){\makebox(0,0)[lb]{\smash{{\SetFigFont{11}{13.2}{\rmdefault}{\bfdefault}{\updefault}{\color[rgb]{0,0,0}$e_1$}%
}}}}
\put(6076,239){\makebox(0,0)[lb]{\smash{{\SetFigFont{11}{13.2}{\rmdefault}{\bfdefault}{\updefault}{\color[rgb]{0,0,0}$s_1$}%
}}}}
\put(10276,-1711){\makebox(0,0)[lb]{\smash{{\SetFigFont{11}{13.2}{\rmdefault}{\bfdefault}{\updefault}{\color[rgb]{0,0,0}$e_0$}%
}}}}
\put(9751,-2236){\makebox(0,0)[lb]{\smash{{\SetFigFont{11}{13.2}{\rmdefault}{\bfdefault}{\updefault}{\color[rgb]{0,0,0}$s_5$}%
}}}}
\put(9001,-2386){\makebox(0,0)[lb]{\smash{{\SetFigFont{11}{13.2}{\rmdefault}{\bfdefault}{\updefault}{\color[rgb]{0,0,0}$e_5$}%
}}}}
\put(8251,-2236){\makebox(0,0)[lb]{\smash{{\SetFigFont{11}{13.2}{\rmdefault}{\bfdefault}{\updefault}{\color[rgb]{0,0,0}$s_4$}%
}}}}
\put(7651,-1711){\makebox(0,0)[lb]{\smash{{\SetFigFont{11}{13.2}{\rmdefault}{\bfdefault}{\updefault}{\color[rgb]{0,0,0}$e_4$}%
}}}}
\put(8176,239){\makebox(0,0)[lb]{\smash{{\SetFigFont{11}{13.2}{\rmdefault}{\bfdefault}{\updefault}{\color[rgb]{0,0,0}$s_2$}%
}}}}
\put(676,-2836){\makebox(0,0)[lb]{\smash{{\SetFigFont{14}{16.8}{\rmdefault}{\bfdefault}{\updefault}{\color[rgb]{0,0,0}$[0 \ 1][2[3 \ 4]5]$}%
}}}}
\put(3047,-4250){\makebox(0,0)[lb]{\smash{{\SetFigFont{11}{13.2}{\rmdefault}{\bfdefault}{\updefault}{\color[rgb]{0,0,0}$e_1$}%
}}}}
\put(2521,-3719){\makebox(0,0)[lb]{\smash{{\SetFigFont{11}{13.2}{\rmdefault}{\bfdefault}{\updefault}{\color[rgb]{0,0,0}$s_1$}%
}}}}
\put(1779,-3509){\makebox(0,0)[lb]{\smash{{\SetFigFont{11}{13.2}{\rmdefault}{\bfdefault}{\updefault}{\color[rgb]{0,0,0}$e_2$}%
}}}}
\put(1038,-3721){\makebox(0,0)[lb]{\smash{{\SetFigFont{11}{13.2}{\rmdefault}{\bfdefault}{\updefault}{\color[rgb]{0,0,0}$s_2$}%
}}}}
\put(457,-4250){\makebox(0,0)[lb]{\smash{{\SetFigFont{11}{13.2}{\rmdefault}{\bfdefault}{\updefault}{\color[rgb]{0,0,0}$e_3$}%
}}}}
\put(507,-5572){\makebox(0,0)[lb]{\smash{{\SetFigFont{11}{13.2}{\rmdefault}{\bfdefault}{\updefault}{\color[rgb]{0,0,0}$e_4$}%
}}}}
\put(347,-4937){\makebox(0,0)[lb]{\smash{{\SetFigFont{11}{13.2}{\rmdefault}{\bfdefault}{\updefault}{\color[rgb]{0,0,0}$s_3$}%
}}}}
\put(6866,-4933){\makebox(0,0)[lb]{\smash{{\SetFigFont{11}{13.2}{\rmdefault}{\bfdefault}{\updefault}{\color[rgb]{0,0,0}$s_0$}%
}}}}
\put(5380,-3510){\makebox(0,0)[lb]{\smash{{\SetFigFont{11}{13.2}{\rmdefault}{\bfdefault}{\updefault}{\color[rgb]{0,0,0}$e_2$}%
}}}}
\put(4056,-4250){\makebox(0,0)[lb]{\smash{{\SetFigFont{11}{13.2}{\rmdefault}{\bfdefault}{\updefault}{\color[rgb]{0,0,0}$e_3$}%
}}}}
\put(4107,-5572){\makebox(0,0)[lb]{\smash{{\SetFigFont{11}{13.2}{\rmdefault}{\bfdefault}{\updefault}{\color[rgb]{0,0,0}$e_4$}%
}}}}
\put(5377,-6314){\makebox(0,0)[lb]{\smash{{\SetFigFont{11}{13.2}{\rmdefault}{\bfdefault}{\updefault}{\color[rgb]{0,0,0}$e_5$}%
}}}}
\put(6597,-5572){\makebox(0,0)[lb]{\smash{{\SetFigFont{11}{13.2}{\rmdefault}{\bfdefault}{\updefault}{\color[rgb]{0,0,0}$e_0$}%
}}}}
\put(3948,-4937){\makebox(0,0)[lb]{\smash{{\SetFigFont{11}{13.2}{\rmdefault}{\bfdefault}{\updefault}{\color[rgb]{0,0,0}$s_3$}%
}}}}
\put(6178,-6096){\makebox(0,0)[lb]{\smash{{\SetFigFont{11}{13.2}{\rmdefault}{\bfdefault}{\updefault}{\color[rgb]{0,0,0}$s_5$}%
}}}}
\put(3226,-4936){\makebox(0,0)[lb]{\smash{{\SetFigFont{11}{13.2}{\rmdefault}{\bfdefault}{\updefault}{\color[rgb]{0,0,0}$s_0$}%
}}}}
\put(1051,-6211){\makebox(0,0)[lb]{\smash{{\SetFigFont{11}{13.2}{\rmdefault}{\bfdefault}{\updefault}{\color[rgb]{0,0,0}$s_4$}%
}}}}
\put(1801,-6361){\makebox(0,0)[lb]{\smash{{\SetFigFont{11}{13.2}{\rmdefault}{\bfdefault}{\updefault}{\color[rgb]{0,0,0}$e_5$}%
}}}}
\put(2551,-6136){\makebox(0,0)[lb]{\smash{{\SetFigFont{11}{13.2}{\rmdefault}{\bfdefault}{\updefault}{\color[rgb]{0,0,0}$s_5$}%
}}}}
\put(3076,-5686){\makebox(0,0)[lb]{\smash{{\SetFigFont{11}{13.2}{\rmdefault}{\bfdefault}{\updefault}{\color[rgb]{0,0,0}$e_0$}%
}}}}
\put(6676,-4261){\makebox(0,0)[lb]{\smash{{\SetFigFont{11}{13.2}{\rmdefault}{\bfdefault}{\updefault}{\color[rgb]{0,0,0}$e_1$}%
}}}}
\put(6151,-3661){\makebox(0,0)[lb]{\smash{{\SetFigFont{11}{13.2}{\rmdefault}{\bfdefault}{\updefault}{\color[rgb]{0,0,0}$s_1$}%
}}}}
\put(4651,-3661){\makebox(0,0)[lb]{\smash{{\SetFigFont{11}{13.2}{\rmdefault}{\bfdefault}{\updefault}{\color[rgb]{0,0,0}$s_2$}%
}}}}
\put(4651,-6136){\makebox(0,0)[lb]{\smash{{\SetFigFont{11}{13.2}{\rmdefault}{\bfdefault}{\updefault}{\color[rgb]{0,0,0}$s_4$}%
}}}}
\put(4501,-2836){\makebox(0,0)[lb]{\smash{{\SetFigFont{14}{16.8}{\rmdefault}{\bfdefault}{\updefault}{\color[rgb]{0,0,0}$[0[1 \ 2]3][4 \ 5]$}%
}}}}
\put(8101,-2836){\makebox(0,0)[lb]{\smash{{\SetFigFont{14}{16.8}{\rmdefault}{\bfdefault}{\updefault}{\color[rgb]{0,0,0}$[0[1[2 \ 3]4]5]$}%
}}}}
\put(4501,-6736){\makebox(0,0)[lb]{\smash{{\SetFigFont{14}{16.8}{\rmdefault}{\bfdefault}{\updefault}{\color[rgb]{0,0,0}$[0 \ 1][2 \ 3][4 \ 5]$}%
}}}}
\put(901,-6736){\makebox(0,0)[lb]{\smash{{\SetFigFont{14}{16.8}{\rmdefault}{\bfdefault}{\updefault}{\color[rgb]{0,0,0}$[0[1 \ 2][3 \ 4]5]$}%
}}}}
\end{picture}%

%% file: d4wsepal.pstex_t
\begin{picture}(0,0)%
\includegraphics{d4wsepal.pstex}%
\end{picture}%
\setlength{\unitlength}{3947sp}%
\begingroup\makeatletter\ifx\SetFigFont\undefined%
\gdef\SetFigFont#1#2#3#4#5{%
  \reset@font\fontsize{#1}{#2pt}%
  \fontfamily{#3}\fontseries{#4}\fontshape{#5}%
  \selectfont}%
\fi\endgroup%
\begin{picture}(6498,3531)(286,-2927)
\put(301,-996){\makebox(0,0)[lb]{\smash{{\SetFigFont{11}{13.2}{\rmdefault}{\bfdefault}{\updefault}{\color[rgb]{0,0,0}$s_3$}%
}}}}
\put(2999,-307){\makebox(0,0)[lb]{\smash{{\SetFigFont{11}{13.2}{\rmdefault}{\bfdefault}{\updefault}{\color[rgb]{0,0,0}$e_1$}%
}}}}
\put(2470,222){\makebox(0,0)[lb]{\smash{{\SetFigFont{11}{13.2}{\rmdefault}{\bfdefault}{\updefault}{\color[rgb]{0,0,0}$s_1$}%
}}}}
\put(1730,433){\makebox(0,0)[lb]{\smash{{\SetFigFont{11}{13.2}{\rmdefault}{\bfdefault}{\updefault}{\color[rgb]{0,0,0}$e_2$}%
}}}}
\put(990,222){\makebox(0,0)[lb]{\smash{{\SetFigFont{11}{13.2}{\rmdefault}{\bfdefault}{\updefault}{\color[rgb]{0,0,0}$s_2$}%
}}}}
\put(407,-307){\makebox(0,0)[lb]{\smash{{\SetFigFont{11}{13.2}{\rmdefault}{\bfdefault}{\updefault}{\color[rgb]{0,0,0}$e_3$}%
}}}}
\put(459,-1630){\makebox(0,0)[lb]{\smash{{\SetFigFont{11}{13.2}{\rmdefault}{\bfdefault}{\updefault}{\color[rgb]{0,0,0}$e_4$}%
}}}}
\put(1730,-2372){\makebox(0,0)[lb]{\smash{{\SetFigFont{11}{13.2}{\rmdefault}{\bfdefault}{\updefault}{\color[rgb]{0,0,0}$e_5$}%
}}}}
\put(2470,-2159){\makebox(0,0)[lb]{\smash{{\SetFigFont{11}{13.2}{\rmdefault}{\bfdefault}{\updefault}{\color[rgb]{0,0,0}$s_5$}%
}}}}
\put(3151,-1036){\makebox(0,0)[lb]{\smash{{\SetFigFont{11}{13.2}{\rmdefault}{\bfdefault}{\updefault}{\color[rgb]{0,0,0}$s_0$}%
}}}}
\put(1051,-2236){\makebox(0,0)[lb]{\smash{{\SetFigFont{11}{13.2}{\rmdefault}{\bfdefault}{\updefault}{\color[rgb]{0,0,0}$s_4$}%
}}}}
\put(3001,-1636){\makebox(0,0)[lb]{\smash{{\SetFigFont{11}{13.2}{\rmdefault}{\bfdefault}{\updefault}{\color[rgb]{0,0,0}$e_0$}%
}}}}
\put(6764,-1020){\makebox(0,0)[lb]{\smash{{\SetFigFont{11}{13.2}{\rmdefault}{\bfdefault}{\updefault}{\color[rgb]{0,0,0}$s_0$}%
}}}}
\put(5281,408){\makebox(0,0)[lb]{\smash{{\SetFigFont{11}{13.2}{\rmdefault}{\bfdefault}{\updefault}{\color[rgb]{0,0,0}$e_2$}%
}}}}
\put(4540,198){\makebox(0,0)[lb]{\smash{{\SetFigFont{11}{13.2}{\rmdefault}{\bfdefault}{\updefault}{\color[rgb]{0,0,0}$s_2$}%
}}}}
\put(4645,-2237){\makebox(0,0)[lb]{\smash{{\SetFigFont{11}{13.2}{\rmdefault}{\bfdefault}{\updefault}{\color[rgb]{0,0,0}$s_4$}%
}}}}
\put(5281,-2395){\makebox(0,0)[lb]{\smash{{\SetFigFont{11}{13.2}{\rmdefault}{\bfdefault}{\updefault}{\color[rgb]{0,0,0}$e_5$}%
}}}}
\put(3853,-1020){\makebox(0,0)[lb]{\smash{{\SetFigFont{11}{13.2}{\rmdefault}{\bfdefault}{\updefault}{\color[rgb]{0,0,0}$s_3$}%
}}}}
\put(4051,-286){\makebox(0,0)[lb]{\smash{{\SetFigFont{11}{13.2}{\rmdefault}{\bfdefault}{\updefault}{\color[rgb]{0,0,0}$e_3$}%
}}}}
\put(4051,-1636){\makebox(0,0)[lb]{\smash{{\SetFigFont{11}{13.2}{\rmdefault}{\bfdefault}{\updefault}{\color[rgb]{0,0,0}$e_4$}%
}}}}
\put(6076,-2236){\makebox(0,0)[lb]{\smash{{\SetFigFont{11}{13.2}{\rmdefault}{\bfdefault}{\updefault}{\color[rgb]{0,0,0}$s_5$}%
}}}}
\put(6601,-1711){\makebox(0,0)[lb]{\smash{{\SetFigFont{11}{13.2}{\rmdefault}{\bfdefault}{\updefault}{\color[rgb]{0,0,0}$e_0$}%
}}}}
\put(6676,-286){\makebox(0,0)[lb]{\smash{{\SetFigFont{11}{13.2}{\rmdefault}{\bfdefault}{\updefault}{\color[rgb]{0,0,0}$e_1$}%
}}}}
\put(6076,239){\makebox(0,0)[lb]{\smash{{\SetFigFont{11}{13.2}{\rmdefault}{\bfdefault}{\updefault}{\color[rgb]{0,0,0}$s_1$}%
}}}}
\put(1051,-2836){\makebox(0,0)[lb]{\smash{{\SetFigFont{14}{16.8}{\rmdefault}{\bfdefault}{\updefault}{\color[rgb]{0,0,0}$[0 \ 1][2 \ 3 \ 4 \ 5]$}%
}}}}
\put(4801,-2836){\makebox(0,0)[lb]{\smash{{\SetFigFont{14}{16.8}{\rmdefault}{\bfdefault}{\updefault}{\color[rgb]{0,0,0}$0[1 \ 2]3[4 \ 5]$}%
}}}}
\end{picture}%

%% file: d4wsepalnhom.pstex_t
\begin{picture}(0,0)%
\includegraphics{d4wsepalnhom.pstex}%
\end{picture}%
\setlength{\unitlength}{3947sp}%
\begingroup\makeatletter\ifx\SetFigFont\undefined%
\gdef\SetFigFont#1#2#3#4#5{%
  \reset@font\fontsize{#1}{#2pt}%
  \fontfamily{#3}\fontseries{#4}\fontshape{#5}%
  \selectfont}%
\fi\endgroup%
\begin{picture}(10754,3686)(61,-2932)
\put(3258,-961){\makebox(0,0)[lb]{\smash{{\SetFigFont{12}{14.4}{\rmdefault}{\bfdefault}{\updefault}{\color[rgb]{0,0,0}$s_0$}%
}}}}
\put(1726,571){\makebox(0,0)[lb]{\smash{{\SetFigFont{12}{14.4}{\rmdefault}{\bfdefault}{\updefault}{\color[rgb]{0,0,0}$e_2$}%
}}}}
\put( 76,-961){\makebox(0,0)[lb]{\smash{{\SetFigFont{12}{14.4}{\rmdefault}{\bfdefault}{\updefault}{\color[rgb]{0,0,0}$s_3$}%
}}}}
\put(901,-2316){\makebox(0,0)[lb]{\smash{{\SetFigFont{12}{14.4}{\rmdefault}{\bfdefault}{\updefault}{\color[rgb]{0,0,0}$s_4$}%
}}}}
\put(1726,-2493){\makebox(0,0)[lb]{\smash{{\SetFigFont{12}{14.4}{\rmdefault}{\bfdefault}{\updefault}{\color[rgb]{0,0,0}$e_5$}%
}}}}
\put(2492,-2316){\makebox(0,0)[lb]{\smash{{\SetFigFont{12}{14.4}{\rmdefault}{\bfdefault}{\updefault}{\color[rgb]{0,0,0}$s_5$}%
}}}}
\put(3140,-1727){\makebox(0,0)[lb]{\smash{{\SetFigFont{12}{14.4}{\rmdefault}{\bfdefault}{\updefault}{\color[rgb]{0,0,0}$e_0$}%
}}}}
\put(253,-1727){\makebox(0,0)[lb]{\smash{{\SetFigFont{12}{14.4}{\rmdefault}{\bfdefault}{\updefault}{\color[rgb]{0,0,0}$e_4$}%
}}}}
\put(3081,-195){\makebox(0,0)[lb]{\smash{{\SetFigFont{12}{14.4}{\rmdefault}{\bfdefault}{\updefault}{\color[rgb]{0,0,0}$e_1$}%
}}}}
\put(194,-195){\makebox(0,0)[lb]{\smash{{\SetFigFont{12}{14.4}{\rmdefault}{\bfdefault}{\updefault}{\color[rgb]{0,0,0}$e_3$}%
}}}}
\put(842,394){\makebox(0,0)[lb]{\smash{{\SetFigFont{12}{14.4}{\rmdefault}{\bfdefault}{\updefault}{\color[rgb]{0,0,0}$s_2$}%
}}}}
\put(2492,394){\makebox(0,0)[lb]{\smash{{\SetFigFont{12}{14.4}{\rmdefault}{\bfdefault}{\updefault}{\color[rgb]{0,0,0}$s_1$}%
}}}}
\put(10800,-961){\makebox(0,0)[lb]{\smash{{\SetFigFont{12}{14.4}{\rmdefault}{\bfdefault}{\updefault}{\color[rgb]{0,0,0}$s_0$}%
}}}}
\put(9268,512){\makebox(0,0)[lb]{\smash{{\SetFigFont{12}{14.4}{\rmdefault}{\bfdefault}{\updefault}{\color[rgb]{0,0,0}$e_2$}%
}}}}
\put(7619,-1020){\makebox(0,0)[lb]{\smash{{\SetFigFont{12}{14.4}{\rmdefault}{\bfdefault}{\updefault}{\color[rgb]{0,0,0}$s_3$}%
}}}}
\put(8443,-2375){\makebox(0,0)[lb]{\smash{{\SetFigFont{12}{14.4}{\rmdefault}{\bfdefault}{\updefault}{\color[rgb]{0,0,0}$s_4$}%
}}}}
\put(9268,-2552){\makebox(0,0)[lb]{\smash{{\SetFigFont{12}{14.4}{\rmdefault}{\bfdefault}{\updefault}{\color[rgb]{0,0,0}$e_5$}%
}}}}
\put(10034,-2375){\makebox(0,0)[lb]{\smash{{\SetFigFont{12}{14.4}{\rmdefault}{\bfdefault}{\updefault}{\color[rgb]{0,0,0}$s_5$}%
}}}}
\put(10683,-1786){\makebox(0,0)[lb]{\smash{{\SetFigFont{12}{14.4}{\rmdefault}{\bfdefault}{\updefault}{\color[rgb]{0,0,0}$e_0$}%
}}}}
\put(7795,-1786){\makebox(0,0)[lb]{\smash{{\SetFigFont{12}{14.4}{\rmdefault}{\bfdefault}{\updefault}{\color[rgb]{0,0,0}$e_4$}%
}}}}
\put(10624,-254){\makebox(0,0)[lb]{\smash{{\SetFigFont{12}{14.4}{\rmdefault}{\bfdefault}{\updefault}{\color[rgb]{0,0,0}$e_1$}%
}}}}
\put(7736,-254){\makebox(0,0)[lb]{\smash{{\SetFigFont{12}{14.4}{\rmdefault}{\bfdefault}{\updefault}{\color[rgb]{0,0,0}$e_3$}%
}}}}
\put(8385,335){\makebox(0,0)[lb]{\smash{{\SetFigFont{12}{14.4}{\rmdefault}{\bfdefault}{\updefault}{\color[rgb]{0,0,0}$s_2$}%
}}}}
\put(10034,335){\makebox(0,0)[lb]{\smash{{\SetFigFont{12}{14.4}{\rmdefault}{\bfdefault}{\updefault}{\color[rgb]{0,0,0}$s_1$}%
}}}}
\put(7088,-961){\makebox(0,0)[lb]{\smash{{\SetFigFont{12}{14.4}{\rmdefault}{\bfdefault}{\updefault}{\color[rgb]{0,0,0}$s_0$}%
}}}}
\put(5556,571){\makebox(0,0)[lb]{\smash{{\SetFigFont{12}{14.4}{\rmdefault}{\bfdefault}{\updefault}{\color[rgb]{0,0,0}$e_2$}%
}}}}
\put(3906,-961){\makebox(0,0)[lb]{\smash{{\SetFigFont{12}{14.4}{\rmdefault}{\bfdefault}{\updefault}{\color[rgb]{0,0,0}$s_3$}%
}}}}
\put(4731,-2316){\makebox(0,0)[lb]{\smash{{\SetFigFont{12}{14.4}{\rmdefault}{\bfdefault}{\updefault}{\color[rgb]{0,0,0}$s_4$}%
}}}}
\put(5556,-2493){\makebox(0,0)[lb]{\smash{{\SetFigFont{12}{14.4}{\rmdefault}{\bfdefault}{\updefault}{\color[rgb]{0,0,0}$e_5$}%
}}}}
\put(6322,-2316){\makebox(0,0)[lb]{\smash{{\SetFigFont{12}{14.4}{\rmdefault}{\bfdefault}{\updefault}{\color[rgb]{0,0,0}$s_5$}%
}}}}
\put(6970,-1727){\makebox(0,0)[lb]{\smash{{\SetFigFont{12}{14.4}{\rmdefault}{\bfdefault}{\updefault}{\color[rgb]{0,0,0}$e_0$}%
}}}}
\put(4083,-1727){\makebox(0,0)[lb]{\smash{{\SetFigFont{12}{14.4}{\rmdefault}{\bfdefault}{\updefault}{\color[rgb]{0,0,0}$e_4$}%
}}}}
\put(6911,-195){\makebox(0,0)[lb]{\smash{{\SetFigFont{12}{14.4}{\rmdefault}{\bfdefault}{\updefault}{\color[rgb]{0,0,0}$e_1$}%
}}}}
\put(4024,-195){\makebox(0,0)[lb]{\smash{{\SetFigFont{12}{14.4}{\rmdefault}{\bfdefault}{\updefault}{\color[rgb]{0,0,0}$e_3$}%
}}}}
\put(4672,394){\makebox(0,0)[lb]{\smash{{\SetFigFont{12}{14.4}{\rmdefault}{\bfdefault}{\updefault}{\color[rgb]{0,0,0}$s_2$}%
}}}}
\put(6322,394){\makebox(0,0)[lb]{\smash{{\SetFigFont{12}{14.4}{\rmdefault}{\bfdefault}{\updefault}{\color[rgb]{0,0,0}$s_1$}%
}}}}
\put(1051,-2836){\makebox(0,0)[lb]{\smash{{\SetFigFont{14}{16.8}{\rmdefault}{\bfdefault}{\updefault}{\color[rgb]{0,0,0}$[0(1 \ 2)3][4 \ 5]$}%
}}}}
\put(8701,-2836){\makebox(0,0)[lb]{\smash{{\SetFigFont{14}{16.8}{\rmdefault}{\bfdefault}{\updefault}{\color[rgb]{0,0,0}$(0 \ 1 \ 2 \ 3)[4 \ 5]$}%
}}}}
\put(4951,-2836){\makebox(0,0)[lb]{\smash{{\SetFigFont{14}{16.8}{\rmdefault}{\bfdefault}{\updefault}{\color[rgb]{0,0,0}$(0[1(2 \ 3)4]5)$}%
}}}}
\end{picture}%